# Non-Hydrostatic Pressure Shallow Flows: GPU Implementation Using Finite Volume and Finite Difference Scheme.

C. Escalante *[1], T. Morales de Luna[2], and M.J. Castro[1]

[1]Departamento de Análisis Matemático, Estadística e Investigación Operativa, y Matemática Aplicada, Universidad de Málaga, Spain

[2]Departamento de Matemáticas, Universidad de Córdoba, Spain

June 27, 2018

**Abstract**

We consider the depth-integrated non-hydrostatic system derived by Yamazaki et al. An efficient formally second-order well-balanced hybrid finite volume finite difference numerical scheme is proposed. The scheme consists of a two-step algorithm based on a projection-correction type scheme initially introduced by Chorin-Temam [15]. First, the hyperbolic part of the system is discretized using a *Polynomial Viscosity Matrix* path-conservative finite volume method. Second, the dispersive terms are solved by means of compact finite differences. A new methodology is also presented to handle wave breaking over complex bathymetries. This adapts well to GPU-architectures and guidelines about its GPU implementation are introduced. The method has been applied to idealized and challenging experimental test cases, which shows the efficiency and accuracy of the method.

*Email address: escalante@uma.es; Corresponding author

*Keywords:*
Non-hydrostatic Shallow-Water, Finite Difference, Finite Volume, GPU, Tsunami Simulation, Wave Breaking

# 1 Introduction

When modelling and simulating geophysical flows, the Nonlinear Shallow-Water equations, hereinafter SWE, are often a good choice as an approximation of the Navier-Stokes equations. Nevertheless, SWE do not take into account effects associated with dispersive waves. In recent years, effort has been made in the derivation of relatively simple mathematical models for shallow water flows that include long nonlinear water waves. As computational power increases, Boussinesq-type Models ([1], [5], [20], [29], [32], [33], [21], [39], [40]) become more accessible. This means that one can use more sophisticated models in order to improve the description of reality, despite the higher computational cost.

Moreover, in order to improve nonlinear dispersive properties of the model, information on the vertical structure of the flow should be included. The Boussinesq-type wave equations have prevailed due to their computational efficiency. The main idea is to include non-hydrostatic effects due to the vertical acceleration of the fluid in the depth-averaging process of the equations. For instance, one can assume that both non-linearity and frequency dispersion are weak and of the same order of magnitude. Since the early works of Peregrine [33], several improved and enhanced Boussinesq models have been proposed over the years: Madsen and Sørensen [29], Nwogu [32], Serre Green-Naghdi equations [20], and nonlinear and non-hydrostatic higher order Shallow-Water type models [7], [41] among many others.

One may use different approaches to improve nonlinear dispersive properties of the models: to consider a Taylor expansion of the velocity potential in powers of the vertical coordinate and in terms of the depth-averaged velocity [29] or the particle velocity components $(u, w)$ at a chosen level [32]; to use a better flow resolution in the vertical direction with a multi-layer approach [26]; to include non-hydrostatic effects in the depth-averaging process ([41], [7]).

The development of non-hydrostatic models for coastal water waves has been the topic of many studies over the past 15 years. Non-hydrostatic models are capable of solving many relevant features of coastal water waves, such as dispersion, non-linearity, shoaling, refraction, diffraction, and run-up (see [7, 41, 14, 36]).

The approach used by Yamazaki in [41] has the advantage of including such non-hydrostatic effects while not adding excessive complexity to the model. This is an advantage from the practical point of view and we will use this technique in this paper.

In this work an already proposed non-hydrostatic pressure system is revisited and some new features are proposed that contributes to the development of an accurate and efficient tool for the simulation of dispersive water waves involving breaking waves, wet-dry fronts and propagation of solitary waves over big domains.

In order to deal with wet-dry fronts in an accurate manner, we have found in several works, that the state of the art when dealing with non-hydrostatic pressure systems, consists in to set to zero the non-hydrostatic pressure for some threshold value (see [41]). In this work, thanks to the rewriting of the incompressibility condition, the non-hydrostatic pressure tends automatically to zero when the water height tends to zero. This is a nice feature, since the only treatment in presence of wet-dry fronts is the redefinition of the hydrostatic pressure term for emerging topographies (to avoid non-physical spurious forces) and the use of a desingularization formula for the computation of velocities when small values of water heights appears.

As is well known, in general non-hydrostatic pressure, and in particular the one studied in this paper, can not deal when breaking waves arises. In such situations one must use a breaking mechanism in order to dissipate the amount of energy associated to turbulence effects when breaking. As we will discuss in the present paper, there are two strategies when dealing with it: the first one is to set the non-hydrostatic pressure to be zero, when a breaking wave is detected. This raises the problem that the convergence of the numerical solution is not ensured when the mesh is refined and a global and costly criteria must be considered (see [23]). The second strategy is to introduce a new physical viscosity term to the horizontal momentum equation. This new term introduces a new parabolic term that must be discretized conveniently. In this work we propose a new writing of a classical breaking term, which allows us to solve the final system in an efficient way.

In this work we will also propose a numerical algorithm that is massively parallelizable, and we will implement it on GPUs architectures. The proposed implementation, which is described in the paper, allows us to compute numerical solutions in big computational domains. This is done using a solely graphic card and reaching a speed-up 110 times faster when compared with a sequential code. The proposed numerical scheme and its implementation for the non-hydrostatic pressure system presents efficient computational times, which are notably similar to an efficient implementation of an hydrostatic shallow water code. This can be stated in fact as the main scientific contribution of this work, which is an advancement and improvement in the field of numerical modelling and numerical simulation of dispersive water waves and, in particular, for non-hydrostatic pressure shallow water system.

The paper is organized as follows. In Section 2, the model is described. In Section 3 breaking mechanism is discussed. The reader should keep in mind that detailed small-scale breaking driven physics are not described by the model. This means that one has to include some breaking mechanism in the depth-integrated equations as it is done by an ad-hoc submodel similar to [35]. In Section 4 a numerical scheme is introduced based on a two-step algorithm. On the first step we solve the SWE in conservative form and on the second step we include the non-hydrostatic effects. The extension of the scheme to the 2D case is introduced in Section 5. In Section 6, guidelines for the GPU implementation of the numerical scheme presented in the previous section are given. Finally, in Section 7, some numerical tests including comparisons with laboratory data are shown.

# 2 Governing equations

In [41] a 2D non-hydrostatic model was presented. The governing equations are derived from the incompressible Navier-Stokes equations. The equations are obtained by a process of depth averaging on the vertical direction $z$. Unlike it is done for SWE, the pressure is not assumed hydrostatic. Following Stelling and Zijlema [36] and Casulli [14], total pressure is decomposed into a sum of hydrostatic and non-hydrostatic pressures. In order to provide the dynamic free-surface boundary condition, non-hydrostatic pressure is assumed to be zero at free surface level.

In the process of depth averaging, the vertical velocity is supposed to have linear vertical profile. Moreover, in the vertical momentum equation, the vertical advective and dissipative terms, which are assumed to be small compared with their horizontal counterparts, are neglected.

The resulting $x$, $y$ and $z$ momentum equations as well as the continuity equation described in [41] are

$$
\begin{cases}
h_t + \nabla \cdot \boldsymbol{q} = 0, \\[2ex]
\boldsymbol{q}_t + \operatorname{div}\left(\dfrac{\boldsymbol{q} \otimes \boldsymbol{q}}{h}\right) + \nabla\left(\dfrac{1}{2}gh^2 + \dfrac{1}{2}hp\right) = (gh + p)\,\nabla H - \tau, \\[2ex]
hw_t = p, \\[2ex]
\nabla \cdot \boldsymbol{u} + \dfrac{W_\eta - W_b}{h} = 0,
\end{cases}
\tag{1}
$$

where $t$ is time and $g$ is gravitational acceleration. $\mathbf{u} = (u, v)$ contains the depth averaged velocities components in the $x$ and $y$ directions respectively. $w$ is the depth averaged velocity component in the $z$ direction. $\mathbf{q} = h\mathbf{u}$ is the discharge vector in the $x$ and $y$ directions. $W_\eta$ and $W_b$ are the vertical velocities at the free-surface and bottom. $p$ is the non-hydrostatic pressure at the bottom. The flow depth is $h = \eta + H$ where $\eta$ is the surface elevation measured from the still-water level, $H$ is the still water depth (see Figure 1). Here we use a Manning friction law given by

$$
\tau = gh\frac{n^2 u|u|}{h^{4/3}},
$$

where $n$ is the Gauckler-Manning coefficient (see [30]).

Operators $\nabla$ and div denote the gradient vector field and the divergence respectively in the $(x, y)$ direction. The vertical velocity at the bottom is evaluated from the boundary condition

$$
W_b = -\boldsymbol{u} \cdot \nabla H. \tag{2}
$$

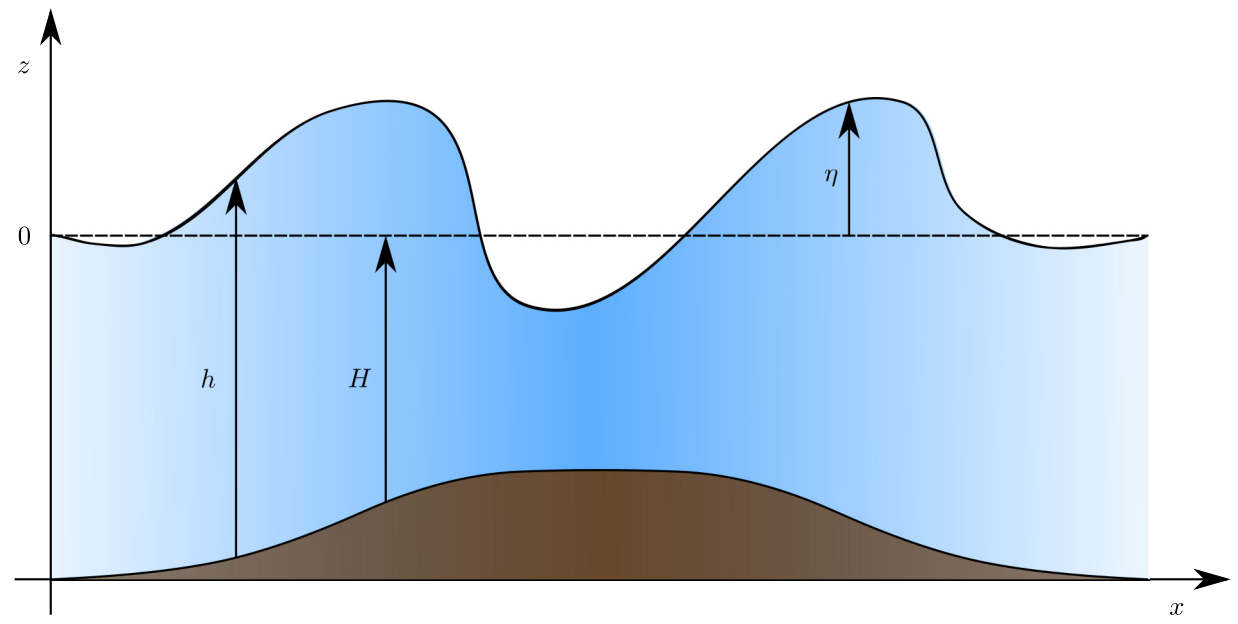

Figure 1: Sketch of the domain for the fluid problem

We will rewrite the system in order to express it in terms of the conserved quantities $h, \boldsymbol{q}$ and $w$. Due to the assumption of a linear profile of the vertical velocity $W$, one has

$$W = W_b + (z + H)\,\frac{W_\eta - W_b}{h},$$

and thus, the integration of $W$ over $z \in [-H, \eta]$ gives

$$w = \frac{1}{h}\int_{-H}^{\eta} W\ dz = \frac{W_\eta + W_b}{2},$$

and therefore,

$$W_\eta = 2w - W_b. \tag{3}$$

Due to the boundary conditions (2) and (3), it holds

$$\frac{W_\eta - W_b}{h} = \frac{w + \boldsymbol{u} \cdot \nabla H}{h/2}, \tag{4}$$

and thus, the last equation in system (1) becomes

$$\nabla \cdot \boldsymbol{u} + \frac{w + \boldsymbol{u} \cdot \nabla H}{h/2} = 0. \tag{5}$$

Finally, equation (5) is multiplied by $h^2$ so that it is rewritten in the form

$$h\nabla \cdot \boldsymbol{q} - \boldsymbol{q} \cdot \nabla\,(2\eta - h) + 2hw = 0, \tag{6}$$

and the system (1) is rewritten as:

$$\begin{cases} h_t + \nabla \cdot \boldsymbol{q} = 0, \\ \boldsymbol{q}_t + \operatorname{div}\left(\dfrac{\boldsymbol{q} \otimes \boldsymbol{q}}{h}\right) + \nabla\left(\dfrac{1}{2}gh^2 + \dfrac{1}{2}hp\right) = (gh + p)\,\nabla H - \tau, \\ hw_t = p, \\ h\nabla \cdot \boldsymbol{q} - \boldsymbol{q} \cdot \nabla\,(2\eta - h) + 2hw = 0. \end{cases} \tag{7}$$

If we consider in system (1) the vertical velocity equation

$$(hw)_t + (hwu)_x = p,$$

then system (7) matches with the one proposed in [7]. In this case, the system verifies an exact energy balance. This property can not be guaranteed for the approach used by Yamazaki in [41], but it has the advantage of not adding excessive complexity to the model. Nevertheless, the numerical scheme proposed in this work can be easily extended to the model proposed in [7]. From the numerical point of view, the results considered in this work do not present relevant differences when comparing both alternatives (see [3, 7])

# 3 Breaking wave modelling

As pointed in [35], in shallow water, complex events can be observed related to turbulent processes. One of these processes corresponds to the breaking of waves near the coast. As it will be seen in the numerical tests proposed in this work, the model presented here cannot describe this process without an additional term which allows the model to dissipate the required amount of energy on such situations. When breaking processes occur, mostly close to shallow areas, two different approaches are usually employed when dispersive Boussinesq-type models are considered.

Close to the coast where breaking starts, the SWE propagates breaking bores at the correct speed, since $kH$ is small, and dissipation of the breaking wave is also well reproduced. Due to that, the simplest way to deal with breaking waves, when considering dispersive systems, consists in neglecting the dispersive part of the equation. This means to force the non-hydrostatic pressure to be zero where breaking occurs. This technique has the advantage that only a breaking criteria is needed to detect this. However, the main disadvantage is that the grid-convergence is not ensured when the mesh is refined, and a global and costly breaking criteria should be taken into account (see [23]).

The other strategy, that will be adopted in this work, consists in dissipation of breaking bores with a diffusive term. Again, a breaking criteria to switch on/off the dissipation is needed. Usually, an eddy viscosity approach (see [35]) solves the matter, where an empirical parameter is defined, based on a quasi-heuristic strategy to determine when the breaking occurs. The main difficulty that presents this mechanism is that usually the diffusive term must be discretized implicitly due to the high order derivatives from the diffusion. Otherwise, it will lead to a severe restriction on the CFL number. As a consequence, an extra linear system has to be solved, losing in efficiency.

We will overcome this challenge in this work.

For the sake of clarity, we will describe the breaking mechanism for the case of one dimensional problems. Let us consider a simple eddy viscosity approach similar to the one introduced in [35], by adding a diffusive term in the horizontal momentum equation of system (7):

$$\begin{cases} h_t + \nabla \cdot \boldsymbol{q} = 0, \\ \boldsymbol{q}_t + \mathrm{div}\left(\dfrac{\boldsymbol{q} \otimes \boldsymbol{q}}{h}\right) + \nabla\left(\dfrac{1}{2}gh^2 + \dfrac{1}{2}hp\right) = (gh+p)\,\nabla H - \tau \\ \qquad\qquad\qquad\qquad\qquad\qquad\qquad\qquad + (\nu h u_x)_x\,, \\ hw_t = p, \\ h\nabla \cdot \boldsymbol{q} - \boldsymbol{q} \cdot \nabla\,(2\eta - h) + 2hw = 0. \end{cases} \tag{8}$$

$\nu$ being the eddy viscosity

$$\nu = Bh|q_x|, \quad B = 1 - \frac{q_x}{U_1}, \tag{9}$$

where

$$U_1 = B_1\sqrt{gh},\ U_2 = B_2\sqrt{gh},$$

denote the flow speeds at the onset and termination of the wave-breaking process and $B_1,\ B_2$ are calibration coefficients that should be fixed through laboratory experiments (see [35]). Wave energy dissipation associated with breaking begins when $|q_x| \geq U_1$ and continues as long as $|q_x| \geq U_2$. The proposed definition of the viscosity $\nu$ requires a positive value of $B$. In order to satisfies that, for negative values of $B$, the viscosity $\nu$ is set to zero.

It is a known fact that using a explicit scheme for a parabolic equation requires a time step restriction of type $\Delta t = \mathcal{O}(\Delta x^2)$. The breaking mechanism has this nature and this would mean a too restrictive time step. This is the reason for choosing an implicit discretization of this term. This can be solved by considering an implicit discretization of the eddy viscosity term, evaluating $\left(\nu^n h_i^{n+1} u_{i,x}^{n+1}\right)_x$ at the right hand side of the momentum discrete equation in (28). The implicit discretization involves solving an extra tridiagonal linear system, leading to a loss of efficiency.

In this work we present, to the best of our knowledge, a new efficient treatment of the eddy viscosity term for depth averaged non-hydrostatic models. To do that, let us rewrite the horizontal momentum equation in (8) as

$$\boldsymbol{q}_t + \operatorname{div}\left(\frac{\boldsymbol{q}\otimes\boldsymbol{q}}{h}\right) + \nabla\left(\frac{1}{2}gh^2 + \frac{1}{2}hp - \nu h u_x\right) = (gh+p)\,\nabla H - \tau, \tag{10}$$

and define

$$p = \widetilde{p} + 2\nu u_x. \tag{11}$$

Thus, replacing $p$ by $\widetilde{p} + 2\nu u_x$, the system (8) can be rewritten as

$$\begin{cases} h_t + \nabla\cdot\boldsymbol{q} = 0, \\ \boldsymbol{q}_t + \operatorname{div}\left(\dfrac{\boldsymbol{q}\otimes\boldsymbol{q}}{h}\right) + \nabla\left(\dfrac{1}{2}gh^2 + \dfrac{1}{2}h\widetilde{p}\right) = (gh+\widetilde{p})\,\nabla H - \tau \\ \qquad\qquad\qquad\qquad\qquad\qquad\qquad\qquad + 2\nu u_x \nabla H, \\ hw_t = \widetilde{p} + 2\nu u_x, \\ h\nabla\cdot\boldsymbol{q} - \boldsymbol{q}\cdot\nabla\left(2\eta - h\right) + 2hw = 0. \end{cases} \tag{12}$$

Note that terms $2\nu u_x \nabla H$, in the horizontal momentum equation, and $2\nu u_x$, in the vertical velocity equation, are essentially first order derivatives of $u$, and can be discretized explicitly without the aforementioned severe restriction on the CFL condition. That gives us an efficient discretization of the eddy viscosity terms.

**Remark 1** *Reinterpretation of the eddy viscosity approach:*

- *Let us consider the vertical component of the stress-tensor*

  $$\tau_{zz} = 2\widetilde{\nu}\partial_z W,$$

  *where $\widetilde{\nu}(x,z,t)$ is a positive function. Now, we use the same process carried out in [41] to depth-average the vertical momentum equation. To do so, let us integrate the vertical component of the stress-tensor along $z \in [-H, \eta]$:*

  $$\int_{-H}^{\eta} \partial_z \tau_{zz}\; dz = 2\int_{-H}^{\eta} \partial_z\widetilde{\nu}\partial_z W + \widetilde{\nu}\partial_{zz} W\; dz.$$

*Since we assume a linear vertical profile for the vertical velocity $W$, then $\partial_{zz}W = 0$ and $\partial_z W$ does not depend on $z$ and thus*

$$\int_{-H}^{\eta} \partial_z \tau_{zz}\ dz = 2\partial_z W \int_{-H}^{\eta} \partial_z \widetilde{\nu}\ dz.$$

*Using again the linearity of the vertical profile for $W$, we get $\partial_z W = \dfrac{W_\eta - W_b}{h}$. From equation (4) and using the last equation in (1) we have that $\partial_z W = -u_x$. Thus,*

$$\int_{-H}^{\eta} \partial_z \tau_{zz}\ dz = -2u_x \int_{-H}^{\eta} \partial_z \widetilde{\nu}\ dz. \tag{13}$$

*Finally, it remains to choose a closure for $\int_{-H}^{\eta} \partial_z \widetilde{\nu}\ dz$ in the system with the described depth-averaged vertical component of the stress-tensor.*

- *If we choose in (12)*

$$\nu = -\int_{-H}^{\eta} \partial_z \widetilde{\nu}\ dz,$$

*then we get the term $2\nu u_x$ introduced in the vertical momentum equation in (12).*

# 4 Numerical scheme

System (8), in the one-dimensional case, can be written in the compact form

$$\begin{cases} \boldsymbol{U}_t + \left(\boldsymbol{F}_{SW}(\boldsymbol{U})\right)_x - \boldsymbol{G}_{SW}(\boldsymbol{U})H_x = \boldsymbol{\mathcal{T}}_{NH}(h, h_x, H, H_x, p, p_x) - \boldsymbol{\tau}, \\ \qquad\qquad\qquad\qquad\qquad\qquad\quad + \boldsymbol{\mathcal{R}}_{\boldsymbol{u}}(\boldsymbol{U}, \boldsymbol{U}_x, H_x) \\ hw_t = p + \mathcal{R}_w(\boldsymbol{U}, \boldsymbol{U}_x), \\ \mathcal{B}(\boldsymbol{U}, \boldsymbol{U}_x, H, H_x, w) = 0, \end{cases} \tag{14}$$

where we introduce the notation

$$\boldsymbol{U} = \begin{pmatrix} h \\ q \end{pmatrix},\ \boldsymbol{F}_{SW}(\boldsymbol{U}) = \begin{pmatrix} q \\ \dfrac{q^2}{h} + \dfrac{1}{2}gh^2 \end{pmatrix},\ \boldsymbol{G}_{SW}(\boldsymbol{U}) = \begin{pmatrix} 0 \\ gh \end{pmatrix},$$

$$\boldsymbol{\mathcal{T}}_{NH}(h, h_x, H, H_x, p, p_x) = \begin{pmatrix} 0 \\ -\dfrac{1}{2}\left(hp_x + p(2\eta - h)_x\right) \end{pmatrix},$$

Finally,

$$\mathcal{B}(\boldsymbol{U}, \boldsymbol{U}_x, H, H_x, w) = hq_x - q\left(2\eta - h\right)_x + 2hw,$$

where

$$\boldsymbol{U}_x = \begin{pmatrix} h_x \\ q_x \end{pmatrix}.$$

and the friction and breaking terms are given by

$$\boldsymbol{\tau} = \begin{pmatrix} 0 \\ \tau \end{pmatrix}, \boldsymbol{\mathcal{R}}_{\boldsymbol{u}}(\boldsymbol{U}, \boldsymbol{U}_x, H_x) = \begin{pmatrix} 0 \\ 2\nu u_x H_x \end{pmatrix}, \mathcal{R}_w(\boldsymbol{U}, \boldsymbol{U}_x) = 2\nu u_x,$$

where $\nu$ is defined by (9).

We describe now the numerical scheme used to discretize the 1D system (14). To do so, we will use a projection method based on the idea introduced in [15]. We shall solve first the hyperbolic problem (SWE). Then, in a second step, non-hydrostatic terms will be taken into account.

The SWE written in vector conservative form is given by

$$\boldsymbol{U}_t + \left(\boldsymbol{F}_{SW}(\boldsymbol{U})\right)_x = \boldsymbol{G}_{SW}(\boldsymbol{U})H_x. \tag{15}$$

The system is solved numerically by using a finite volume method. In particular, an efficient second-order well-balanced *Polynomial Viscosity Matrix* (PVM) path-conservative finite volume method [8] is applied. As usual, we consider a set of finite volume cells $I_i = [x_{i-1/2}, x_{i+1/2}]$ with lengths $\Delta x_i$ and define

$$\boldsymbol{U}_i(t) = \frac{1}{\Delta x_i}\int_{I_i} \mathbf{U}(x,t)dx,$$

the cell average of the function $\boldsymbol{U}(x,t)$ on cell $I_i$ at time $t$.

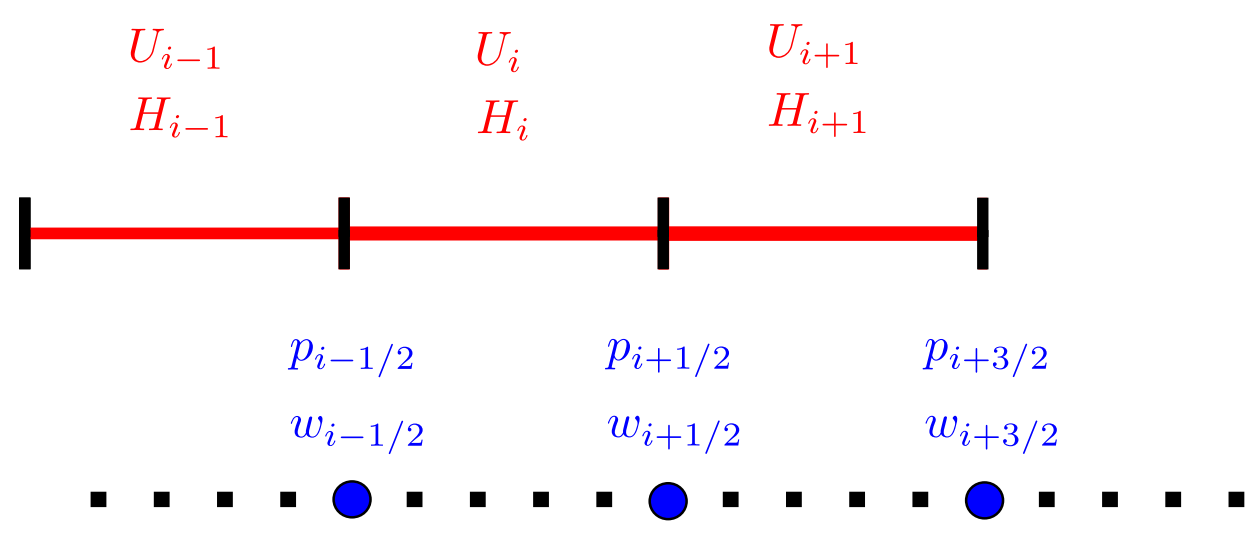

Figure 2: Numerical scheme stencil. Up: finite volume mesh. Down: staggered mesh for finite differences.

Regarding non-hydrostatic terms, we consider a staggered-grid (see Figure 2) formed by the points $x_{i-1/2}$, $x_{i+1/2}$ of the interfaces for each cell $I_i$, and denote the point values of the functions $p$ and $w$ on point $x_{i+1/2}$ at time $t$ by

$$p_{i+1/2}(t) = p(x_{i+1/2}, t), \;\; w_{i+1/2}(t) = w(x_{i+1/2}, t).$$

Non-hydrostatic terms will be approximated by second order compact finite differences. The resulting ODE system is discretized using a *Total Variation Diminishing* (TVD) Runge-Kutta method [19]. For the sake of clarity, only a first order discretization in time will be described. The source terms corresponding to friction terms are discretized semi-implicitly. The breaking terms are discretized explicitly using second order finite differences. Thus, friction terms are neglected and only flux, and source terms are considered.

## 4.1 Finite volume discretization for the SWE

For the sake of simplicity we shall consider a constant cell length $\Delta x$. A first order path-conservative PVM scheme for system (14) reads as follows (see [8]):

$$\mathbf{U}_i'(t) = -\frac{1}{\Delta x}\left(D^-_{i+1/2}(t) + D^+_{i-1/2}(t)\right), \tag{16}$$

where, avoiding the time dependence,

$$\begin{aligned} D^\pm_{i+1/2} &= D^\pm_{i+1/2}(\mathbf{U}_i, \mathbf{U}_{i+1}, H_i, H_{i+1}) = \\ &= \frac{1}{2}\left(\mathbf{F}(\mathbf{U}_{i+1}) - \mathbf{F}(\mathbf{U}_i) - \mathbf{G}_{i+1/2}\left(H_{i+1} - H_i\right)\right) \\ &\pm \frac{1}{2} Q_{i+1/2}\left(\left(\mathbf{U}_{i+1} - \mathbf{U}_i\right) - A^{-1}_{i+1/2}\mathbf{G}_{i+1/2}\left(H_{i+1} - H_i\right)\right), \end{aligned} \tag{17}$$

where

$$\mathbf{G}_{i+1/2} = \begin{pmatrix} 0 \\ gh_{i+1/2} \end{pmatrix},$$

and

$$A_{i+1/2} = \begin{pmatrix} 0 & 1 \\ -u_{i+1/2}^2 + gh_{i+1/2} & 2u_{i+1/2} \end{pmatrix}$$

is the Roe Matrix associated to the flux $\mathbf{F}(\mathbf{U})$ from the SWE, being

$$h_{i+1/2} = \frac{h_i + h_{i+1}}{2}, \; u_{i+1/2} = \frac{u_i\sqrt{h_i} + u_{i+1}\sqrt{h_{i+1}}}{\sqrt{h_i} + \sqrt{h_{i+1}}}.$$

$Q_{i+1/2}$ is the viscosity matrix associated to the numerical method. For PVM schemes, $Q_{i+1/2}$ is obtained by a polynomial evaluation of the Roe Matrix.
In this work, the viscosity matrix is defined as

$$Q_{i+1/2} = \alpha_0 Id + \alpha_1 A_{i+1/2},$$

being

$$\alpha_0 = \frac{S_R|S_L| - S_L|S_R|}{S_R - S_L}, \quad \alpha_1 = \frac{|S_R| - |S_L|}{S_R - S_L},$$

where

$$S_L = \min\left(u_{i+1/2} - \sqrt{gh_{i+1/2}}, u_i - \sqrt{gh_i}\right),$$

$$S_R = \max\left(u_{i+1/2} + \sqrt{gh_{i+1/2}}, u_{i+1} + \sqrt{gh_{i+1}}\right).$$

Under this choice, $D^{\pm}_{i+1/2}$ read as

$$\begin{aligned} D^{\pm}_{i+1/2} = & \frac{1}{2}\left(\mathbf{F}(\mathbf{U}_{i+1}) - \mathbf{F}(\mathbf{U}_i) - \mathbf{G}_{i+1/2}\left(H_{i+1} - H_i\right)\right) \\ & \pm \frac{1}{2}\left(\alpha_0 Id + \alpha_1 A_{i+1/2}\right)\left((\mathbf{U}_{i+1} - \mathbf{U}_i) - A^{-1}_{i+1/2}\mathbf{G}_{i+1/2}\left(H_{i+1} - H_i\right)\right). \end{aligned}$$

The scheme is a path-conservative extension of the Harten-Lax-van Leer (HLL) scheme ([22])

Note that the above expression is not well defined for the resonant case when $A_{i+1/2}$ is not invertible. This problem can be avoided following the strategy described in [12], where $A_{i+1/2}$ is replaced by

$$A^*_{i+1/2} = \begin{pmatrix} 0 & 1 \\ gh_{i+1/2} & 0 \end{pmatrix}.$$

This choice will be made in general. On the one hand, this makes the scheme simpler. On the other hand, it avoids singularities at critical points. This means that there is no need to check whether we are near a critical point or not. The counterpart is that the scheme is only well-balanced for lake at rest steady states.

Using this particular choice, the numerical scheme reads as

$$D^{\pm}_{i+1/2} = \frac{1}{2}\left((1 \pm \alpha_1)R_{i+1/2} \pm \alpha_0\left(U_{i+1} - U_i - (H_{i+1} - H_i)\,\boldsymbol{e}_1\right)\right), \tag{18}$$

where

$$R_{i+1/2} = F_c(\mathbf{U}_{i+1}) - F_c(\mathbf{U}_i) + T_{p,i+1/2},\ \boldsymbol{e}_1 = \begin{pmatrix} 1 \\ 0 \end{pmatrix},$$

being

$$F_c(\mathbf{U}_i) = \begin{pmatrix} q_i \\ \dfrac{q_i^2}{h_i} \end{pmatrix}, T_{p,i+1/2} = \begin{pmatrix} 0 \\ gh_{i+1/2}(\eta_{i+1} - \eta_i) \end{pmatrix} \tag{19}$$

the corresponding discretization of convective and pressure terms for the SWE.

Second order in space is obtained following [11] by combining a MUSCL reconstruction operator (see [25]) with the PVM scheme presented above, that can be written as

$$\mathbf{U}_i'(t) = -\frac{1}{\Delta x}\left(D^-_{i+1/2}(t) + D^+_{i-1/2}(t) + \mathcal{I}_i(t)\right), \tag{20}$$

where

$$D^{\pm}_{i+1/2} = D^{\pm}_{i+1/2}(\mathbf{U}^-_{i+1/2}, \mathbf{U}^+_{i+1/2}, H^-_{i+1/2}, H^+_{i+1/2}), \tag{21}$$

$$\mathcal{I}_i = \boldsymbol{F}(\boldsymbol{U}^-_{i+1/2}) - \boldsymbol{F}(\boldsymbol{U}^+_{i-1/2}) - \boldsymbol{G}(\boldsymbol{U}_i)\left(H^-_{i+1/2} - H^+_{i-1/2}\right). \tag{22}$$

The vector $\boldsymbol{U}^{\pm}_{i+1/2}$ is defined by the reconstructed variables $h^{\pm}_{i+1/2}$, $\eta^{\pm}_{i+1/2}$, $u^{\pm}_{i+1/2}$, to the left ($^-$) and right ($^+$) of the inter-cell $x_{i+1/2}$, from the cell averages applying a MUSCL reconstruction operator (see [25]), combined with a minmod limiter. The MUSCL reconstruction operator takes into account the positivity of the water height. Finally, the variable $H^{\pm}_{i+1/2}$ is recovered from $H^{\pm}_{i+1/2} = h^{\pm}_{i+1/2} - \eta^{\pm}_{i+1/2}$. This procedure allows the scheme to be well-balanced for the water at rest solutions, and to deal with emerging topographies: since the variable $\eta$ is reconstructed instead oh $H$, in these situations, the gradient of $\eta$ is set to zero, to avoid spurious non-physical pressure forces (see [10]).

**Remark 2** *Concerning the well-balancing properties, the numerical scheme considered in this work (first or second order) is well-balanced for the water at rest solution and is linearly $L^\infty$-stable under the usual CFL condition, that is*

$$\Delta t < CFL \frac{\Delta x}{|\lambda_{max}|},\ 0 < CFL \leq 1,\ |\lambda_{max}| = \max_{i\in\{1,\dots,N\}} \left\{ |u_i| + \sqrt{g h_i} \right\}. \quad (23)$$

## 4.2 Finite difference discretization for the non-hydrostatic terms

In this Subsection, non-hydrostatic variables $p$, $w$ will be discretized using second order compact finite differences. In order to obtain point value approximations for the non-hydrostatic variables $p_{i+1/2},\ w_{i+1/2}$, and skipping notation in time, the operator $\mathcal{B}(\boldsymbol{U}, \boldsymbol{U}_x, H, H_x, w)$ will be approximated for every point $x_{i+1/2}$ of the staggered grid (Figure 2) by

$$\begin{aligned} &\mathcal{B}(\boldsymbol{U}_{i+1/2}, \boldsymbol{U}_{x,i+1/2}, H_{i+1/2}, H_{x,i+1/2}, w_{i+1/2}) = \\ &h_{i+1/2} q_{x,i+1/2} - q_{i+1/2}\left(2\eta_{x,i+1/2} - h_{x,i+1/2}\right) + 2h_{i+1/2} w_{i+1/2}, \end{aligned} \quad (24)$$

where we will use second order point value approximations of $\boldsymbol{U}, \boldsymbol{U}_x, H$ and $H_x$, on the staggered-grid. They will be computed from the approximations of the average values on the cell $I_i,\ I_{i+1}$ as follows:

$$\begin{aligned} h_{i+1/2} = \frac{h_{i+1} + h_i}{2},\ h_{x,i+1/2} = \frac{h_{i+1} - h_i}{\Delta x},\ \eta_{x,i+1/2} = \frac{\eta_{i+1} - \eta_i}{\Delta x}, \\ q_{i+1/2} = \frac{q_{i+1} + q_i}{2},\ q_{x,i+1/2} = \frac{q_{i+1} - q_i}{\Delta x}. \end{aligned} \quad (25)$$

In a similar way, a second order point value approximation in the center of the cell will be used for $\boldsymbol{\mathcal{T}}_{NH}$, computed as

$$\boldsymbol{\mathcal{T}}_{NH}(h_i, h_{x,i}, H_i, H_{x,i}, p_i, p_{x,i}) = \begin{pmatrix} 0 \\ -\frac{1}{2}\left(h_i p_{x,i} + p_i(2\eta_{x,i} - h_{x,i})\right) \end{pmatrix}. \quad (26)$$

Here

$$h_{x,i} = \frac{h_{i+1} - h_{i-1}}{2\Delta x}, \;\; \eta_{x,i} = \frac{\eta_{i+1} - \eta_{i-1}}{2\Delta x},$$
$$p_i = \frac{p_{i+1/2} + p_{i-1/2}}{2}, \;\; p_{x,i} = \frac{p_{i+1/2} - p_{i-1/2}}{\Delta x}, \tag{27}$$

are second order point value approximations in the middle of the cell $I_i$, which are a second order approximation of the averaged variables.

## 4.3 Final numerical scheme

Assume given time steps $\Delta t^n$, and denote $t^n = \sum_{k\leq n} \Delta t^k$ and $\mathbf{U}_i(t^n) = \mathbf{U}_i^n$, $p_{i+1/2}(t^n) = p_{i+1/2}^n$, $w_{i+1/2}(t^n) = w_{i+1/2}^n$. The numerical scheme proposed can be summarized as follows:
In a first stage, SWE approximations are solved. Let us define $U_i^{n+1/2}$ as the averaged values of $U$ on cell $I_i$ at time $t^n$ for the SWE as detailed in the Subsection 4.1.
In a second stage, we consider the system

$$\begin{cases} \boldsymbol{U}_i^{n+1} = \boldsymbol{U}_i^{n+1/2} + \Delta t \boldsymbol{\mathcal{R}_u}(\boldsymbol{U}^{n+1/2}, \boldsymbol{U}_{x,i}^{n+1/2}, H_{x,i}) \\ \qquad\qquad + \Delta t \boldsymbol{\mathcal{T}}_{NH}(h_i^{n+1}, h_{x,i}^{n+1}, H_i, H_{x,i}, p_i^{n+1}, p_{x,i}^{n+1}), \\ w_{i+1/2}^{n+1} = w_{i+1/2}^n + \Delta t \mathcal{R}_w(\boldsymbol{U}_i^{n+1/2}, \boldsymbol{U}_{x,i}^{n+1/2}) + \Delta t \dfrac{p_{i+1/2}^{n+1}}{h_{i+1/2}^{n+1}}, \\ \mathcal{B}(\boldsymbol{U}_{i+1/2}^{n+1}, \boldsymbol{U}_{x,i+1/2}^{n+1}, H_{i+1/2}, H_{x,i+1/2}, w_{i+1/2}^{n+1}) = 0, \end{cases} \tag{28}$$

where

$$\mathcal{B}(\boldsymbol{U}_{i+1/2}^{n+1}, \boldsymbol{U}_{x,i+1/2}^{n+1}, H_{i+1/2}, H_{x,i+1/2}, w_{i+1/2}^{n+1})$$

is given by (24) and

$$\boldsymbol{\mathcal{T}}_{NH}(h_i^{n+1}, h_{x,i}^{n+1}, H_i, H_{x,i}, p_i^{n+1}, p_{x,i}^{n+1})$$

is given by (26). Finally $\boldsymbol{U}_{x,i}^{n+1/2}$ and $H_{x,i}$ appearing in the breaking terms are computed as it was done in Subsection 4.2 from the point value approximations in the middle of the cell $I_i$

$$\boldsymbol{U}_{x,i}^{n+1/2} = \frac{\boldsymbol{U}_{i+1}^{n+1/2} - \boldsymbol{U}_{i-1}^{n+1/2}}{2\Delta x}, \;\; H_{x,i} = \frac{(h_{i+1} - \eta_{i+1}) - (h_{i-1} - \eta_{i-1})}{2\Delta x},$$

which are a second order approximation of the averaged variables.

System (28) leads to solve a linear system

$$\boldsymbol{A}^{n+1/2}\boldsymbol{\mathcal{P}}^{n+1} = \boldsymbol{\mathcal{RHS}}^{n+1/2}, \tag{29}$$

where $\boldsymbol{A}^{n+1/2}$ is a tridiagonal matrix. The matrix $\boldsymbol{A}^{n+1/2}$ as well as the *Right Hand Side* vector $\boldsymbol{\mathcal{RHS}}^{n+1/2}$ are given in Appendix A.1. We would also like to stress the dependency of $\boldsymbol{A}$ and $\boldsymbol{\mathcal{RHS}}$ on the variables $h$ and $\eta$ at the time $n+1/2$. $\boldsymbol{\mathcal{P}}^{n+1}$ is a vector containing the non-hydrostatic pressure values at time $n+1$.

The linear system is efficiently solved using the Thomas algorithm. Then the values $q_i^{n+1/2}$ are corrected with $\boldsymbol{\mathcal{T}}_{NH}(h_i^{n+1}, h_{x,i}^{n+1}, H_i, H_{x,i}, p_i^{n+1}, p_{x,i}^{n+1})$.

The scheme presented here is only first order in time. To get a second order in time discretization, we perform a second order TVD Runge-Kutta approach (see [19]). Therefore, the resulting scheme is second order accurate in space and time. Remark that the usual CFL restriction (23) should be considered.

## 4.4 Boundary conditions

In this work, three types of Boundary Conditions (BC) have been considered: periodic, outflow and generating/absorbing BCs.

1. Periodic BCs: Given the domain subdivided into a set of $N$ cells, cell $I_1$ and $I_N$, which are the extremes of the domain, are considered as the same cell, surrounded by the neighbour cells $I_N$ to the left and $I_2$ to the right. In this case, the matrix is no more tridiagonal and a modification of the Thomas algorithm is used.

2. Outflow BCs: homogeneous Neumann conditions are applied on the left and right boundaries. Since we use a second order MUSCL scheme, the usage of one ghost cell $I_0,\ I_{N+1}$ in each boundary is required in order to determine the values of the closest nodes to the boundary. The values of the variables at the ghost cells are extrapolated from the adjacent cells.

   Nevertheless, reflections at the boundaries might modify the numerical solution at the interior of the domain. As in many other works (see [23, 34] among others), this condition is sometimes supplemented with an absorbing BC described bellow.

3. Generating/absorbing BCs: Periodic wave generation as well as absorbing BCs are achieved by using a generation/relaxation zone method similar to the one proposed in [28].

   Generation/absorption of waves is achieved by simply defining a relaxation coefficient $0 \leq m(x) \leq 1$, and a target solution $(\boldsymbol{U}^*, w^*, p^*)$. Given a width $L_{Rel}$ of the relaxation zone on each boundary, we define $k_{Rel}$ as the first natural number that $k_{Rel}\Delta x \geq L_{Rel}$. The solution within the relaxation zone is then redefined to be, $\forall i \in \{1, \dots, k_{Rel}, N - k_{rel}, \dots N\}$ :

$$
\begin{aligned}
\widetilde{\boldsymbol{U}}_i &= m_i \boldsymbol{U}_i + (1 - m_i)\boldsymbol{U}_i^* \\
\widetilde{w}_{i\pm1/2} &= m_{i\pm1/2} w_{i\pm1/2} + (1 - m_{i\pm1/2}) w^*_{i\pm1/2}, \\
\widetilde{p}_{i\pm1/2} &= m_{i\pm1/2} p_{i\pm1/2} + (1 - m_{i\pm1/2}) p^*_{i\pm1/2},
\end{aligned}
$$

   where $m_i$ is defined as

$$
m_i = \sqrt{1 - \left(\frac{d_i}{L_{Rel}}\right)^2}, \; m_{i\pm1/2} = \frac{m_i + m_{i\pm1}}{2},
$$

   where $d_i$ is the distance between the centre of the cells $I_i$ and $I_1$ (respectively $I_i$ and $I_{N-k}$), in the case of $i \in \{1, \dots, k\}$ (respectively $i \in \{N - k, \dots, N\}$).

   For the numerical experiments we set

$$
L \leq L_{Rel} \leq 1.5L,
$$

   $L$ being the typical wavelength of the outgoing wave.

   Absorbing BC is q particular case where $\boldsymbol{U}^* = w^* = p^* = 0$. This will dump all the waves passing through the boundaries.

## 4.5 Wet-dry treatment

For the computation of $\boldsymbol{U}^{n+1/2}$ in the finite volume discretization of the underlying hyperbolic system, a wet-dry treatment adapting the ideas described in [10] is applied. The key of the numerical treatment for wet-dry fronts with emerging bottom topographies relies in:

- The hydrostatic pressure term in (19) at the horizontal velocity equation is modified for emerging bottoms to avoid spurious pressure forces (see [10]).

- To compute velocities appearing in the numerical scheme from the discharges, one has $u = q/h$. This may present difficulties close to dry areas due to small values of $h$, resulting in large round-off errors. We shall compute the velocities analogously as in [24], applying the desingularization formula

$$u = \frac{\sqrt{2}hq}{\sqrt{h^4 + \max(h^4, \delta^4)}},$$

  which gives the exact value of $u$ for $h \geq \delta$, and gives a smooth transition of $u$ to zero when $h$ tends to zero, with no truncation. In this work we set $\delta = 10^{-5}$ for the numerical tests. A more detailed discussion about the desingularization formula can be seen in [24].

In the second step of the numerical scheme, no special treatment is required due to the rewriting of the incompressibility equations, which has been multiplied by $h^2$, and is expressed in terms of discharges. In presence of wet-dry fonts, the non-hydrostatic pressure vanishes and no artificial truncation up to a threshold value is needed. This is shown in Appendix A.2, where an analysis is carried out for the case of

$$h = \delta, \; q = w = 0, \; H = \alpha x.$$

In such situation, we can assert that the linear system that defines the non-hydrostatic pressure at each step, is always invertible. Since the *Right Hand Side* vector of the linear system vanishes, then the only solution for the homogeneous linear system is that the non-hydrostatic pressure vanishes.

# 5 Numerical scheme in two dimensions

We describe the numerical scheme used to discretize the 2D system (7). The computational domain is decomposed into subsets with a simple geometry, called cells or finite volumes. We will use one common arrangement of the variables, known as the Arakawa C-grid (see Figure 3). This is an extension of the procedure used for the 1D case. Variables $p$ and $w$ will be computed at the intersections of the edges:

$$p_{i+1/2,j+1/2}(t) = p(x_{i+1/2}, y_{j+1/2}, t), \;\; w_{i+1/2,j+1/2}(t) = w(x_{i+1/2}, y_{j+1/2}, t).$$

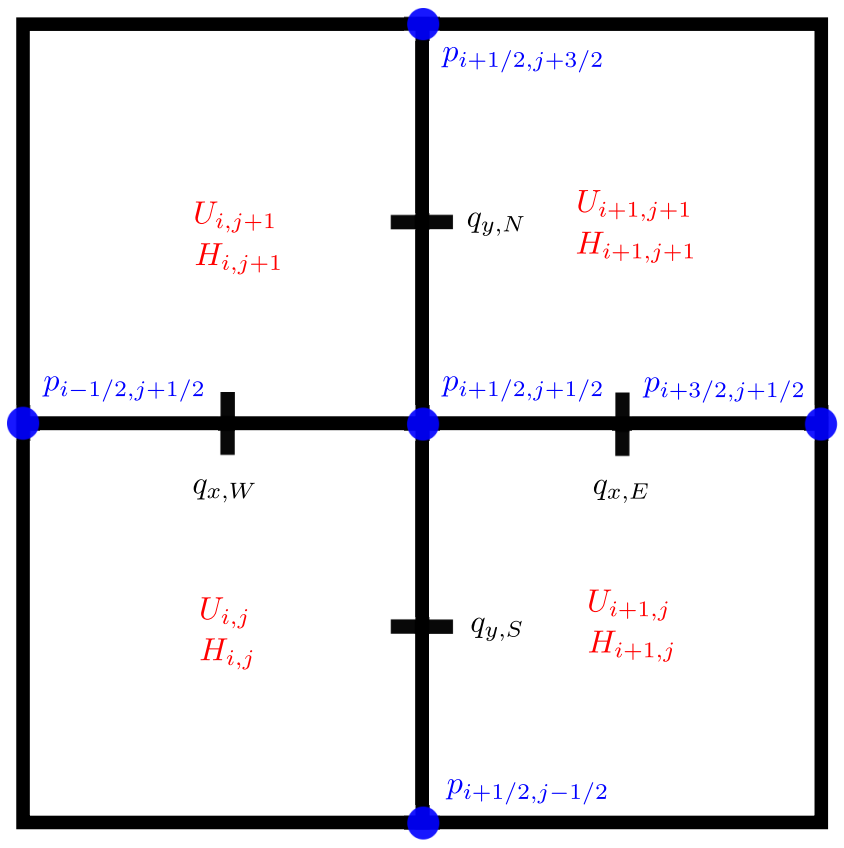

Figure 3: Numerical scheme stencil

As in Section 4, we shall solve first the hyperbolic problem (SWE) and then correct it with the non-hydrostatic terms.

The SWE are solved numerically by using a finite volume method. An efficient second-order well-balanced PVM path-conservative finite volume method is applied following [8]. There, second order in space is obtained following [11] by combining a MUSCL reconstruction operator (see [25]) with the PVM scheme. In particular, we use the 2D extension of the PVM scheme described in Section 3.1 in [18]. We describe here the expression of the second order HLL scheme as:

$$\begin{aligned} \boldsymbol{U}_{i,j}^{n+1/2} &= \boldsymbol{U}_{i,j}^{n} - \frac{1}{|V_{i,j}|} \sum_{k\in\mathcal{N}_{i,j}} |E_{ij(k)}| \mathcal{F}_{ij(k)}^{HLL^-}(\boldsymbol{U}_{i,j}^{n,-}, \boldsymbol{U}_{i,j}^{n,+}, H_{ij}^{-}, H_{ij}^{+}) \\ & \quad - \frac{1}{|V_{i,j}|} \int_{V_{i,j}} \begin{pmatrix} 0 \\ g h_{i,j}^{n}(\boldsymbol{x}) (\eta_x)_{i,j} \\ g h_{i,j}^{n}(\boldsymbol{x}) (\eta_y)_{i,j} \end{pmatrix} \end{aligned} \tag{30}$$

where $\boldsymbol{U}_{i,j}^{n,-}$ and $H_{i,j}^{-}$(respectively $\boldsymbol{U}_{i,j}^{n,+}$ and $H_{i,j}^{+}$) are the values of the reconstruction variables from the cell averages applying a MUSCL reconstruction operator (see [25]) combined with a minmod limiter, at the center of the edge $E_{ij(k)}$ at time $n$, and $(\eta_x)_{i,j}$ (respectively $(\eta_y)_{i,j}$) is the constant approximation of the partial derivative of free surface with respect to $x$ (respectively $y$) at cell $V_{i,j}$ provided by the reconstruction. $h^n_{i,j}(\boldsymbol{x})$ is the reconstruction of the water depth at cell $V_{i,j}$ at time $n$. The integral appearing in (30) is approximated by the mid-point rule.

Non-hydrostatic terms are approximated by second order compact finite differences. The resulting ODE system is discretized using a TVD Runge-Kutta method [19]. The source terms corresponding to friction terms are discretized semi-implicitly. Breaking terms are discretized following the ideas presented in Section 3.

The final numerical scheme is

$$\left\{\begin{array}{l} \boldsymbol{U}^{n+1}_{i,j} = \boldsymbol{U}^{n+1/2}_{i,j} + \Delta t \boldsymbol{\mathcal{T}}_{NH}\left(h^{n+1}, \nabla(h^{n+1}), H, \nabla(H), p^{n+1}, \nabla(p^{n+1})\right)_{i,j}, \\ \\ w^{n+1}_{i+1/2,j+1/2} = w^{n}_{i+1/2,j+1/2} + \Delta t \dfrac{p^{n+1}_{i+1/2,j+1/2}}{h^{n+1}_{i+1/2,j+1/2}}, \\ \\ \mathcal{B}\left(\boldsymbol{U}^{n+1}, \nabla(h^{n+1}), (\nabla\cdot\boldsymbol{Q}^{n+1}), H, \nabla(H), w^{n+1}\right)_{i+1/2,j+1/2} = 0. \end{array}\right. \tag{31}$$

where we denote the vector of the state variables

$$\boldsymbol{U} = \begin{pmatrix} h \\ \boldsymbol{Q} \end{pmatrix}, \quad \boldsymbol{Q} = \begin{pmatrix} q_1 \\ q_2 \end{pmatrix},$$

and $\mathcal{B}$, $\mathcal{T}_{NH}$ defined as in Section 4. $\mathcal{B}$ will be approximated for every point $x_{i+1/2,j+1/2}$ of the staggered-grid. To do that, second order point value approximations of $\widetilde{\boldsymbol{U}}^{n+1}, \nabla(h^{n+1}), (\nabla\cdot\widetilde{\boldsymbol{Q}}^{n+1}), H, \nabla(H)$ and $w^{n+1}$ on the staggered-grid points will be computed from the approximations of the average values on the cell provided in the first SWE finite volume step.

In the same way, a second order point value approximation in the center of the cell will be used for the approximation of $\boldsymbol{\mathcal{T}}_{NH}$.

System (31) leads to solve a penta-diagonal linear system for the unknowns $p^{n+1}_{i+1/2,j+1/2}$. The coefficients of the matrices depend on the variables $h$ and $\eta$, at time $n+1/2$. Since the resulting coefficients of the matrix are too tedious to be given in this paper, we shall omit them. A rigorous analysis of the matrices in general is not an easy task. Nevertheless, in all the numerical computations, we have checked that the matrices are strictly diagonally dominant. Thus, due to the Gershgorin circle theorem, the matrices are non-singular for all the test cases shown in this paper.

The linear system is solved using an iterative Jacobi method combined with a scheduled relaxation method following [2].

Remark that the compactness of the numerical stencil and the easy parallelization of the Jacobi method adapts well to the implementation of the scheme on GPUs architectures. Given $\boldsymbol{P}^{n+1}$ a vector that contains the non-hydrostatic pressure unknowns, to define a convergence criteria, we use

$$\boldsymbol{E}^{n+1,(k+1)} = \|\boldsymbol{P}^{n+1,(k+1)} - \boldsymbol{P}^{n+1,(k)}\|_\infty < \epsilon \tag{32}$$

where $\boldsymbol{P}^{n+1,(k)}$ denotes the $k$-th approximation of $\boldsymbol{P}^{n+1}$ given by the Jacobi algorithm, and $\epsilon$ is a tolerance parameter. It is observed that the Jacobi method converges in a few iterations for the problems tested here.

To get a second order in time discretization, we perform a second order TVD Runge-Kutta approach (see [19]). The details of the scheme can be found in the Appendix A.1.

# 6 GPU implementation

We are mainly interested in the application to real-life problems: simulation in channels, dambreak problems, ocean currents, tsunami propagation, etc. Simulating those phenomena gives place to long time simulations in big computational domains. Thus, extremely efficient implementations are needed to be able to analyze those problems in low computational time.

The numerical scheme presented here exhibits a high potential for data parallelization. This fact suggests the design of parallel implementation of the numerical scheme. NVIDIA has developed the CUDA programming toolkit [31] for modern Graphics Processor Units (GPUs). CUDA includes an extension of the C language and facilitates the programming on GPUs for general purpose applications by preventing the programmer to deal with low level language programming on GPU.

In this section, guidelines for the implementation of the numerical scheme presented in the previous sections are given. The general steps of the parallel implementation are shown in Figure 4. Each step executed on the GPU is assigned to a CUDA *kernel*, which is a function executed on the GPU. Let us describe the main loop of the program. To do so, let us assume that we have at time $t^n$ the values $\boldsymbol{U}_{i,j}^n$ for each volume $V_{i,j}$ and a precomputed stimation $\Delta t^n$. We will also describe the numerical algorithm for the first order in time case.

At the beginning of the algorithm we build the finite volume mesh and the main data structure to be used in GPU. For each volume $V_{i,j}$ we store the state variables in one array of type `double4`[1]. This array contains $h$, $q_1$, $q_2$ and $H$, given by $\boldsymbol{U}_{i,j}^n$. A series of CUDA *kernels* will do the following tasks:

1. **Process fluxes on edges:** In this step, each thread computes the contribution at every edge of two adjacent volumes. This thread will also compute the volume integral appearing in (30) using the mid-point rule. This implementation follows a similar approach to the one applied in [16] and [17]. The edge processing is succesively done in the horizontal and vertical direction, computing even and odd edges separetly. This avoids simultaneous access to the same memory values by two different threads. The computed contributions are stored in an array accumulator of type `double4` with size equal to the number of volumes (see [16] for further details).

   Note that previous computations require the use of the reconstructed values, $\boldsymbol{U}_{i,j}^{n,-}, \boldsymbol{U}_{i,j}^{n,+}$, as well as the reconstructed topography values, $H_{ij}^-, H_{ij}^+$.

2. **Update $\boldsymbol{U}_{i,j}^{n+1/2}$ for each volume**: In this step, each thread will compute the next state $\boldsymbol{U}_{i,j}^{n+1/2}$ for each volume $V_{i,j}$ by using the values stored in accumulator and the precomputed estimation of $\Delta t^n$. Moreover, a local $\Delta t_{ij}^{n+1}$ is computed for each volume from the CFL condition.

3. **Solve the linear system for non-hydrostatic pressure**: In order to solve the linear system (29), we use a Jacobi iterative method. This implementation is matrix-free, as the the coefficients of the matrix are not pre-computed and stored. Instead, the coefficients are computed on the fly, which means less memory usage. For each point of the staggered mesh $(x_{i+1/2}, y_{j+1/2})$, we store the last two iterations of the non-hydrostatic pressure of the Jacobi algorithm, the local error, and the vertical velocity using an array of type `double4`.

[1]The `double4` data type represent structures with four double precision real components

This step is splitted into two parts: first, given $\boldsymbol{P}^{n+1,(k)}$, a *kernel* will perform an iteration of the Jacobi method, obtaining $\boldsymbol{P}^{n+1,(k+1)}$ and $\boldsymbol{E}^{n+1,(k+1)}$. Second, another CUDA *kernel* will compute the minimum of all local errors by applying a reduction algorithm in GPU.

4. **Compute the values $\boldsymbol{U}_{i,j}^{n+1}$ for each volume**: In this step, each thread has access to a given volume and it computes the next state $\boldsymbol{U}_{i,j}^{n+1}$ by using the values of the non-hydrostatic pressure obtained previously.

5. **Get estimation of $\Delta t^{n+1}$**: Similarly to what is done in [16] and [17], the minimum of all the local $\Delta t_{ij}^{n+1}$ values is obtained by applying a reduction algorithm in GPU. This value shall be used as precomputed $\Delta t^{n+1}$ for the next step of the loop.

When considering a second order discretization in time, the steps 1-4 are repeated twice, for each step of the Runge-Kutta method. Finally, the step 5 is done at the end of the temporal evolution for every time step.

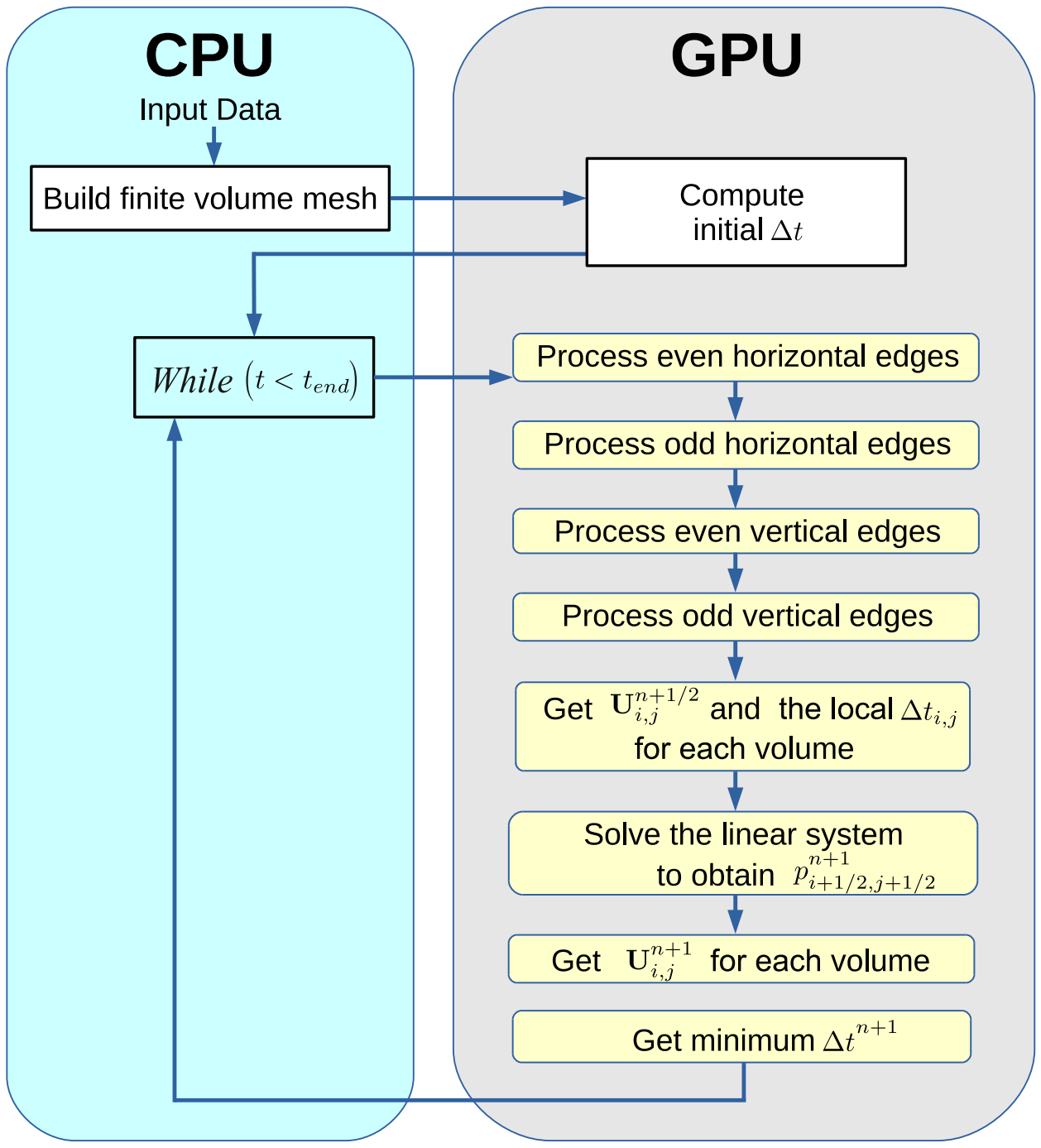

Figure 4: Parallel CUDA implementation.

# 7 Numerical tests and results

## 7.1 Solitary wave propagation in a channel

The propagation of a solitary wave over a long distance is a standard test of the stability and conservative properties of numerical schemes for Boussinesq-type equations ([7], [41], [35], [34], [36], [38]). A solitary wave propagates at constant speed and without change of shape over an horizontal bottom. An approximated expression of a solitary wave for system (7) is given by (see [38])

$$\eta(x,t) = A \cdot Sech^2\left[\sqrt{\frac{3A}{4H^3}}(x-ct))\right],\ u(x,t) = \frac{\sqrt{gH}}{H}\eta(x,t), \tag{33}$$

where $A$ is the amplitude and $c = \sqrt{g(A+H)}$ is the wave propagation velocity.

We simulate the propagation of a solitary wave over a constant depth $H = 1.0\ m$ with $A = 0.1\ m$ in a channel of length 500 $m$ along the $x$ direction. The domain is divided into 5000 cells along the $x$ axis. The final time is 400$s$. We set $CFL = 0.4$ and $g = 1.0\ m/s^2$ . Periodic boundary conditions are considered, and the initial condition is computed using (33).

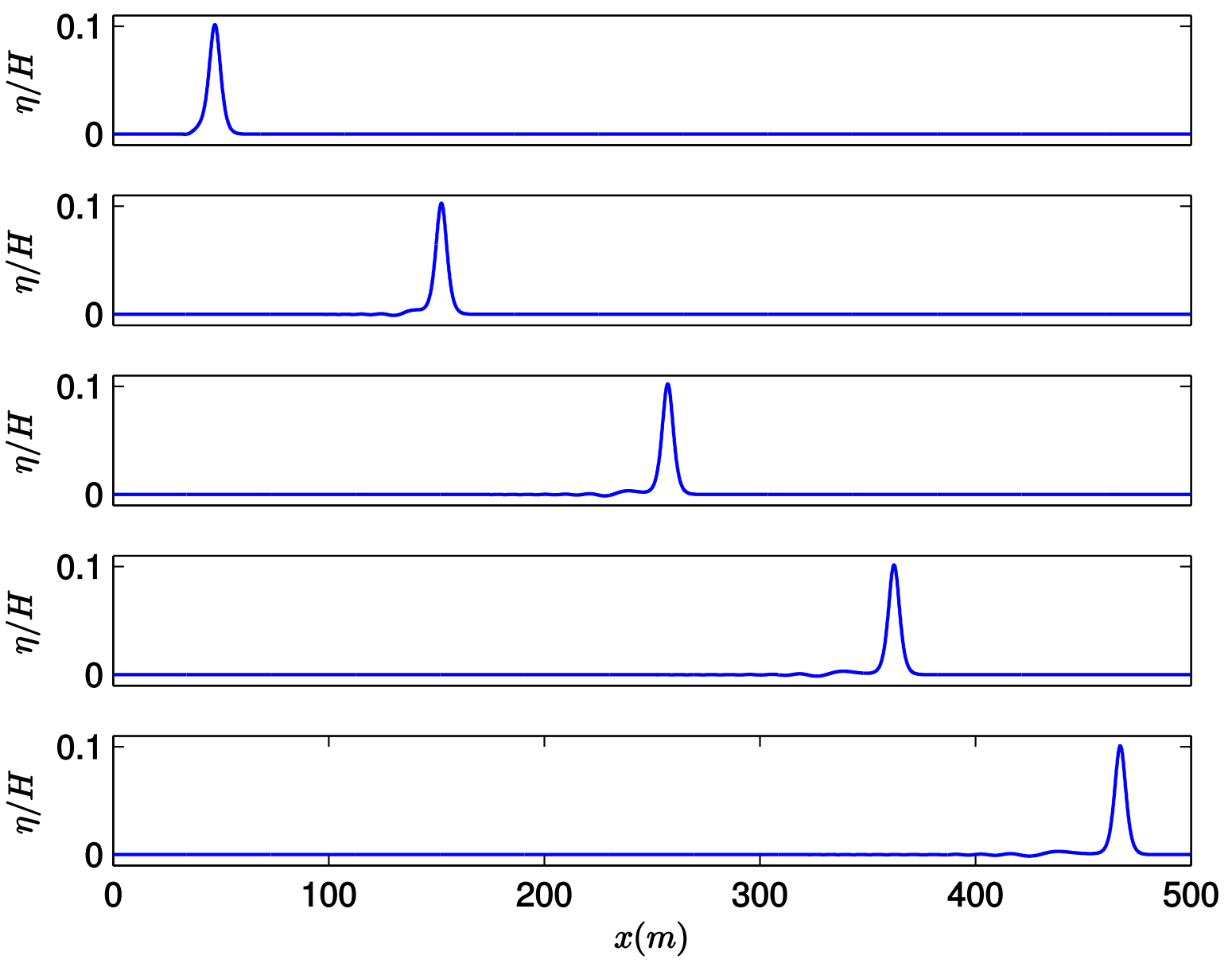

Figure 5: Solitary wave propagation at $T = 0,\ 100,\ 200,\ 300,\ 400\ s$

Figure 5 shows the evolution of the solitary wave at different times. As expected, the wave's shape has not changed and propagates at constant speed (see Figure 6).

| Number of Cells | $L^1$ error $h$ | $L^1$ order $h$ | $L^1$ error $q$ | $L^1$ order $q$ |
|---|---|---|---|---|
| 100 | 2.99E-03 | - | 3.88E-03 | - |
| 200 | 7.19E-04 | 2.06 | 7.44E-04 | 2.02 |
| 400 | 1.78E-04 | 2.01 | 1.78E-04 | 1.98 |
| 800 | 4.51E-05 | 1.98 | 4.23E-05 | 1.96 |
| 1600 | 1.19E-05 | 1.92 | 1.19E-05 | 1.94 |
| 3200 | 3.20E-06 | 1.90 | 3.86E-06 | 1.95 |

Table 1: One-dimensional accuracy test. $L^1$ numerical errors and orders.

Numerical simulations for different grids have been computed up to time $t = 10.0$ $s$ in a channel of length 50 $m$. Table 1 shows the $L^1$ errors and numerical orders of accuracy obtained with $CFL$ number 0.4. Since equation 33 is not an exact solution for system (7), we take as reference solution a numerical simulation at time $t = 10.0$ $s$ for a very fine grid with 12800 cells.

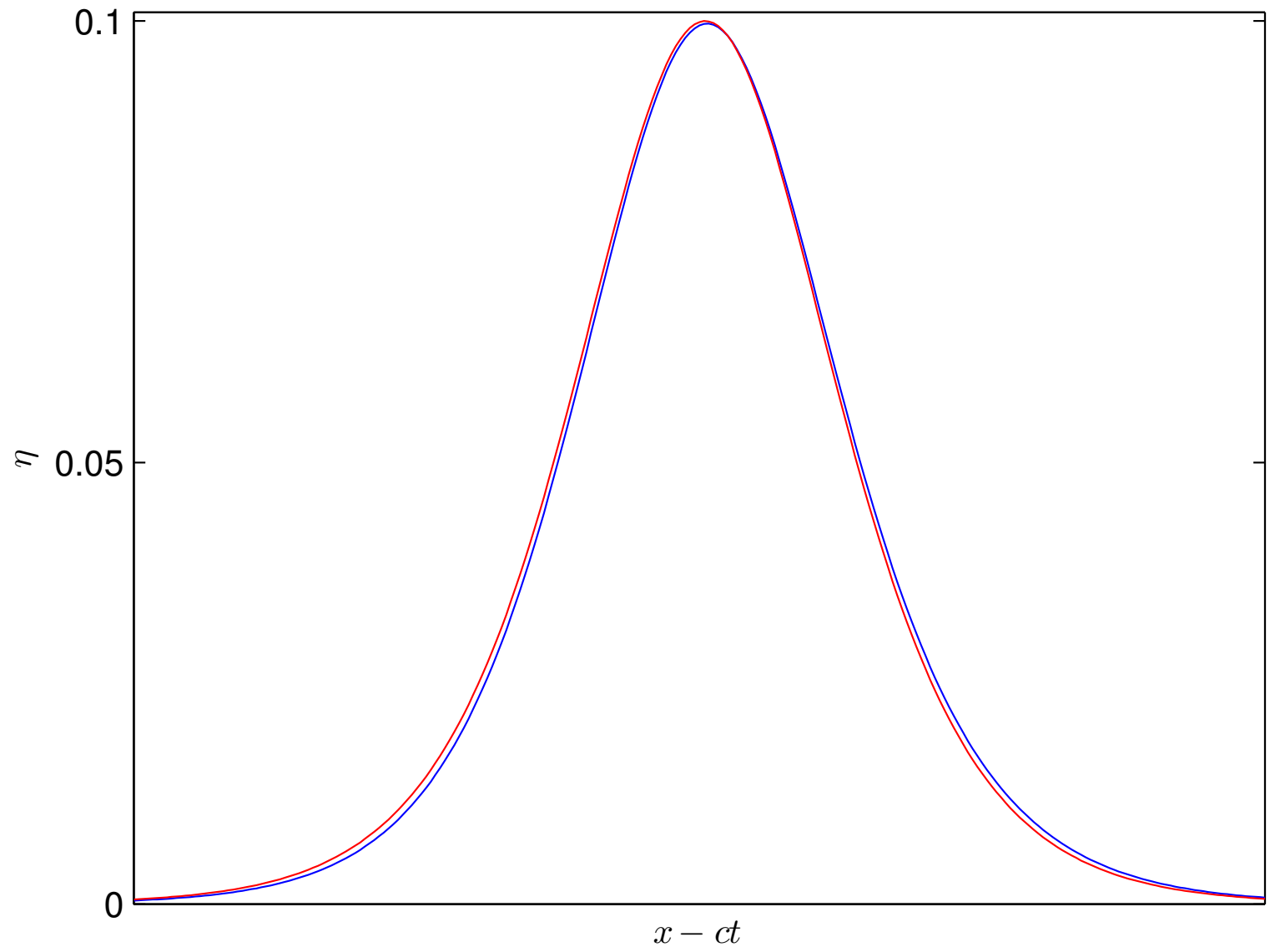

Figure 6: Comparison of analytical (red) and numerical (blue) surface at time $T = 400\ s$

## 7.2 Head-on collision of two solitary waves

The head-on collision of two equal solitary waves is again a common test for the Boussinesq-type models (see [35], [34]). The collision of two solitary waves is equivalent to the reflection of one solitary wave by a vertical wall when viscosity is neglected.

After the interaction of the two waves, one should ideally recover the initial profiles. The collision of the two waves presents additional challenges to the model due to the sudden change of the nonlinear and frequency dispersion characteristics.

We present here the interaction of two solitary waves propagating on a depth of $H = 1\ m$ with amplitude $A = 0.1\ m$. The same computational scenario, same boundary conditions and same expression for the solitary wave (33) as in previous test are taken into account.

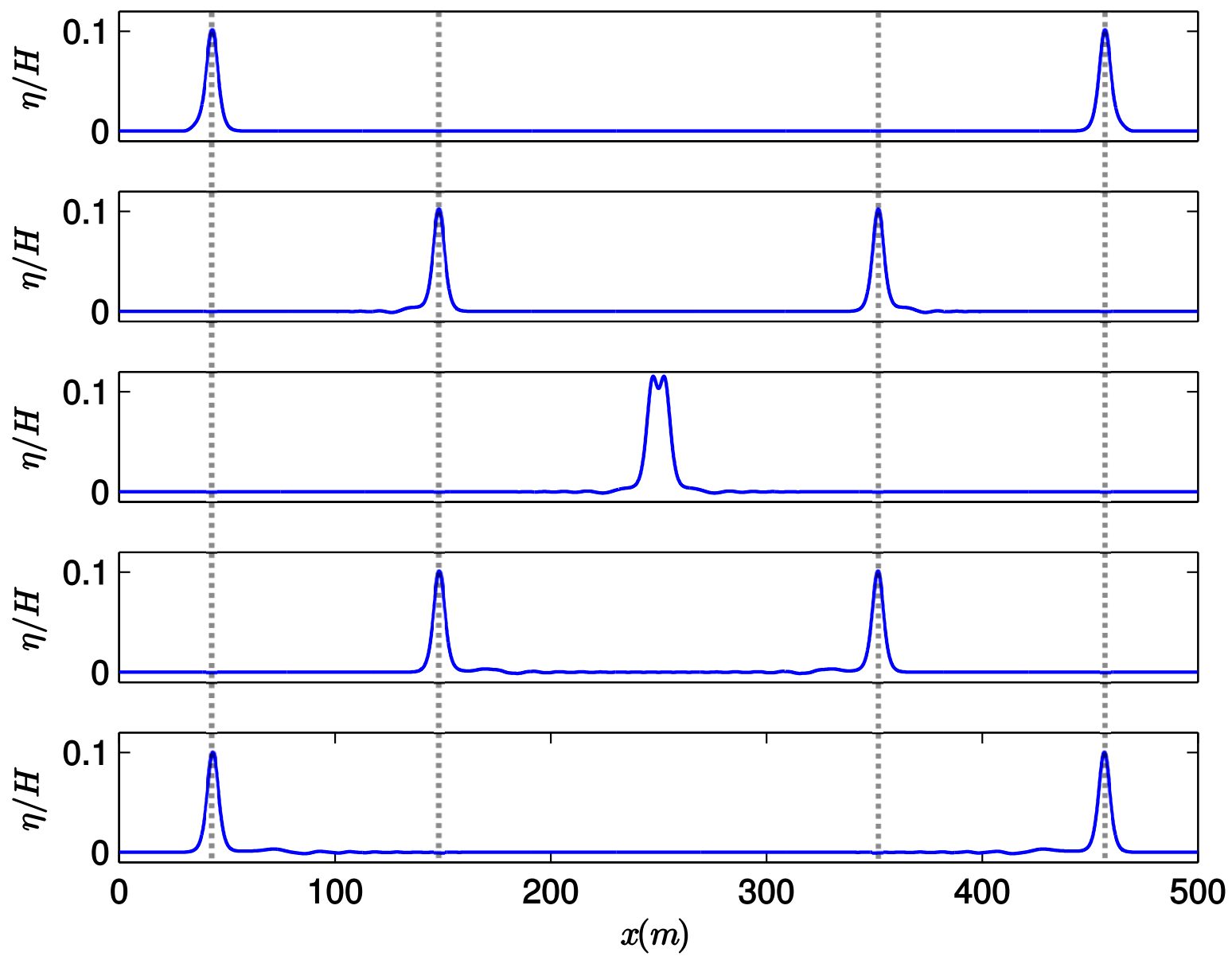

Figure 7: Head-on collision of two solitary waves at $T = 0,\ 100,\ 200,\ 300,\ 400\ s$

Figure 7 shows the collision of the two solitary waves at the midpoint of the domain. After the collision both maintain the initial amplitude and the same speed but in opposite directions.

## 7.3 Periodic waves breaking over a submerged bar

The experiment of plunging breaking periodic waves over a submerged bar by Beji and Battjes [4] is considered here. The numerical test is performed in a one-dimensional channel with a trapezoidal obstacle submerged. Waves in the free surface are measured in seven point stations $S_0, S_1, \ldots, S_6$ ( See Figure 8).

The one-dimensional domain $[0, 25]$ is discretized with $\Delta x = 0.05\ m$. and the bathymetry is defined in the Figure 8.

The velocity $u$ and surface elevation $\eta$ are set initially to 0. The boundary conditions are: free outflow at $x = 25\ m$ and free surface is imposed at $x = 0\ m$ using the data provide by the experiment at $S_0$. The data provided at $S_0$ by the experiment is the free-surface $\eta_{S_0}(t)$ and the velocity $u_{S_0}(t)$. Thus, we use as a target solution for the generating boundary condition (see

Section 4.4)

$$h^*(t) = 0.4 + \eta_{S_0}(t),\ q^*(t) = h^*(t)u_{S_0}(t),\ w^*(t) = 0,\ p^*(t) = 0.$$

The first wave gauge $S_1$ shows that the imposed generating boundary conditions are well implemented, since the match is excellent.

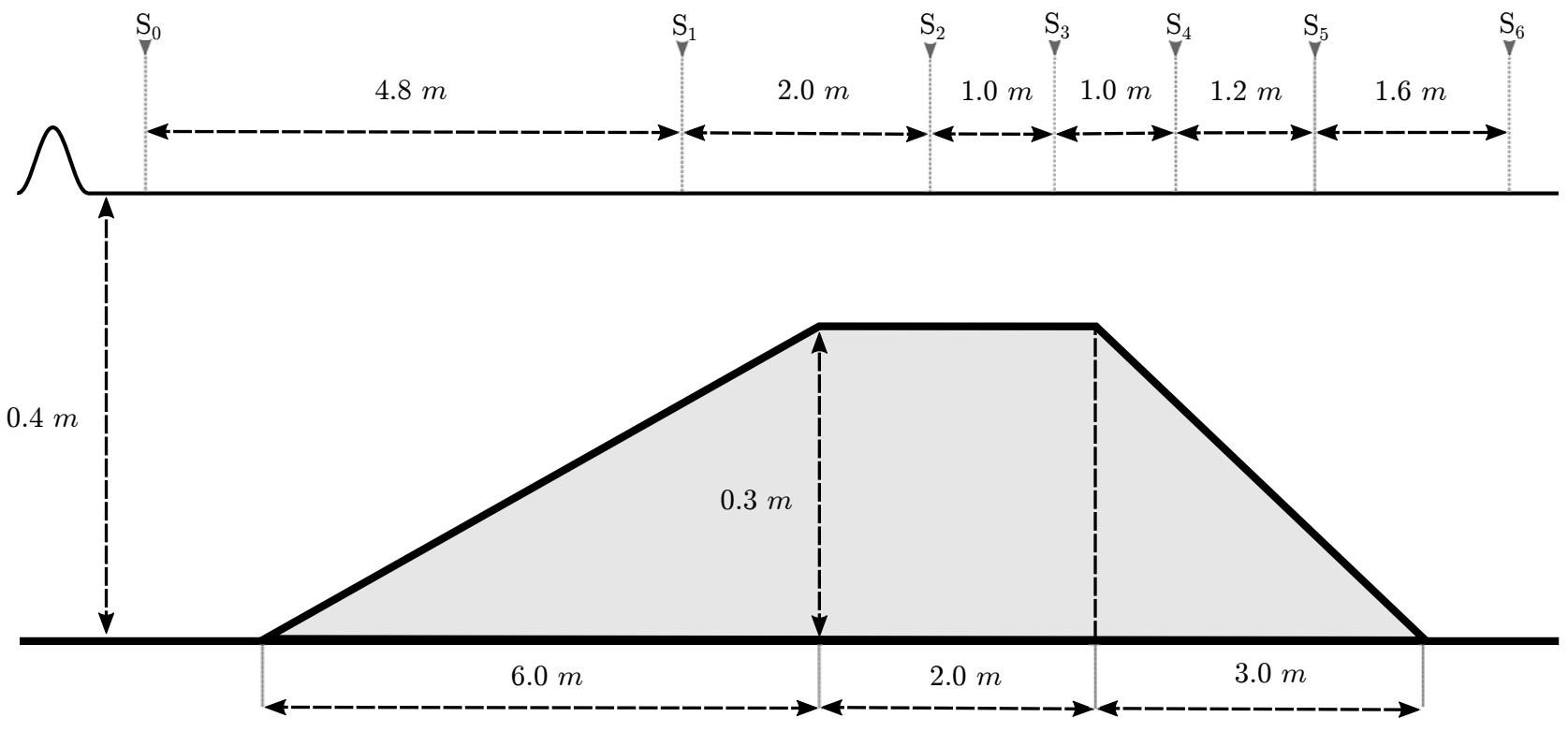

Figure 8: Periodic waves breaking over a submerged bar. Sketch of the topography and layout of the wave gauges

The CFL number is set to 0.9 and $g = 9.81\ m/s^2$. Figure 9 shows the time evolution of the free surface at points $S_1, \ldots, S_6$. The comparison with experimental data emphasizes the need to consider a dispersive model to faithfully capture the shape of the waves near the continental slope. Both amplitude and frequency of the waves are captured on all wave gauges successfully.

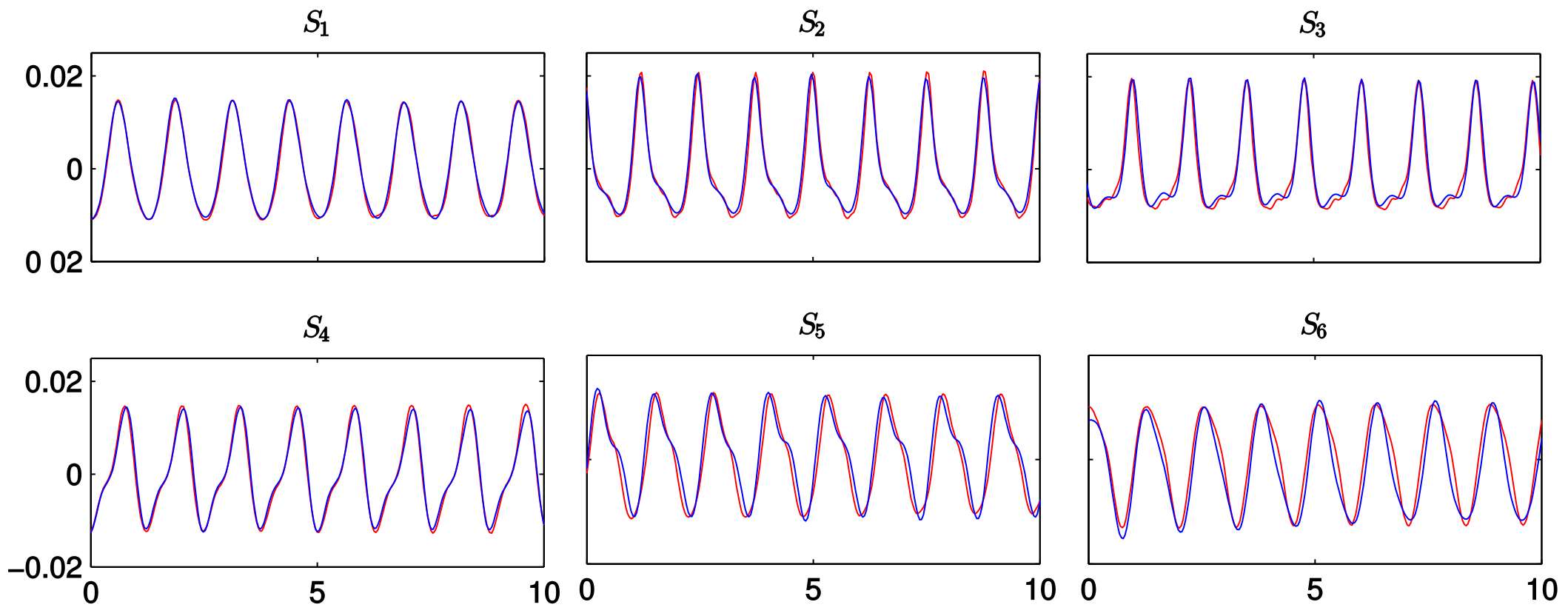

Figure 9: Comparison of data time series (red) and numerical (blue) at wave gauges $S_1,\ S_2,\ S_3,\ S_4,\ S_5,\ S_6$

## 7.4 Solitary wave run-up on a plane beach

Solitary wave run-up on a plane beach is one of the most intensively studied problems in long-wave modeling. Synolakis [37] carried out laboratory experiments for incident solitary waves of multiple relative amplitudes, to study propagation, breaking and run-up over a planar beach with a slope 1 : 19.85. Many researchers have used this data to validate numerical models. With this test case we assess the ability of the model to describe shoreline motions and wave breaking, when it occurs. Experimental data are available in [37] for surface elevation at different times. For this test the still water level is $H = 1\ m$. The bathymetry of the problem is given by Figure 10.

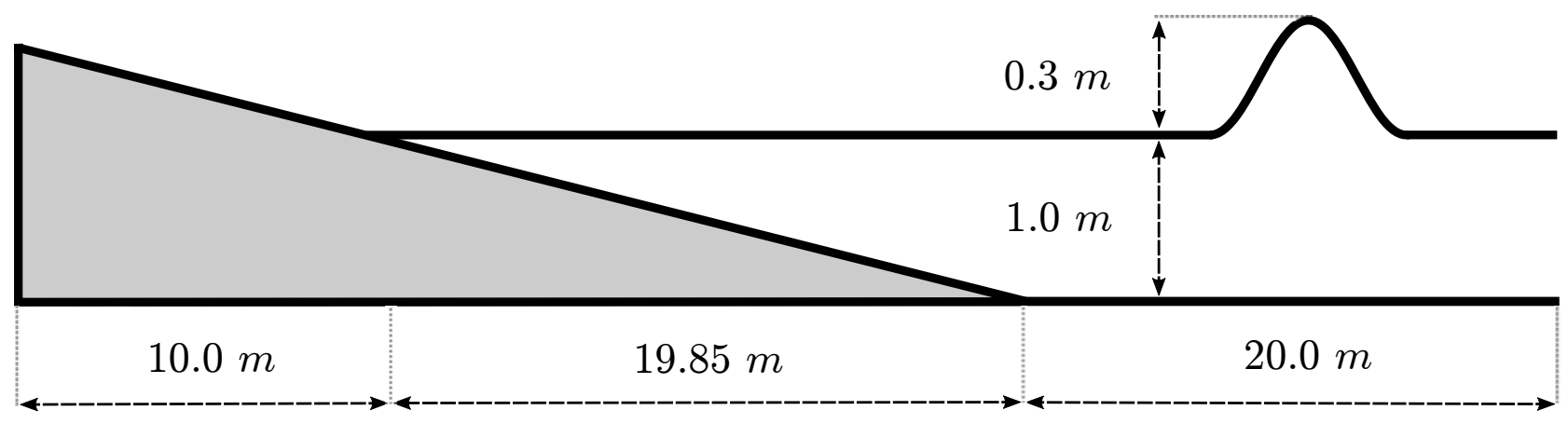

Figure 10: Sketch of the topography

A solitary wave of amplitude $0.3\ m$. is placed at point $x = 25\ m$. given by (33). A Manning coefficient of $n_m = 0.01$ was used in order to define the glass surface roughness used in the experiments. The computational domain is $[-10, 40]$ and the numerical parameters used were $\Delta x = 0.05$ , $CFL = 0.9$ and $g = 9.81\ m/s^2$. Free outflow boundary conditions are imposed.

Figure 11 shows snapshots at different times, $t\sqrt{g/H} = t_0$ where $H = 1$. A good agreement between experimental and simulated data is seen. Here we use the breaking criteria described in Section 3 with $B_1 = 0.15$ and $B_2 = 0.5$. Figure 12 shows the same test case described previously, but this time the breaking mechanism is not considered. In this case, an overshoot value on the amplitude of the wave appears when the mesh is refined. The results are quite satisfactory in favour of the former.

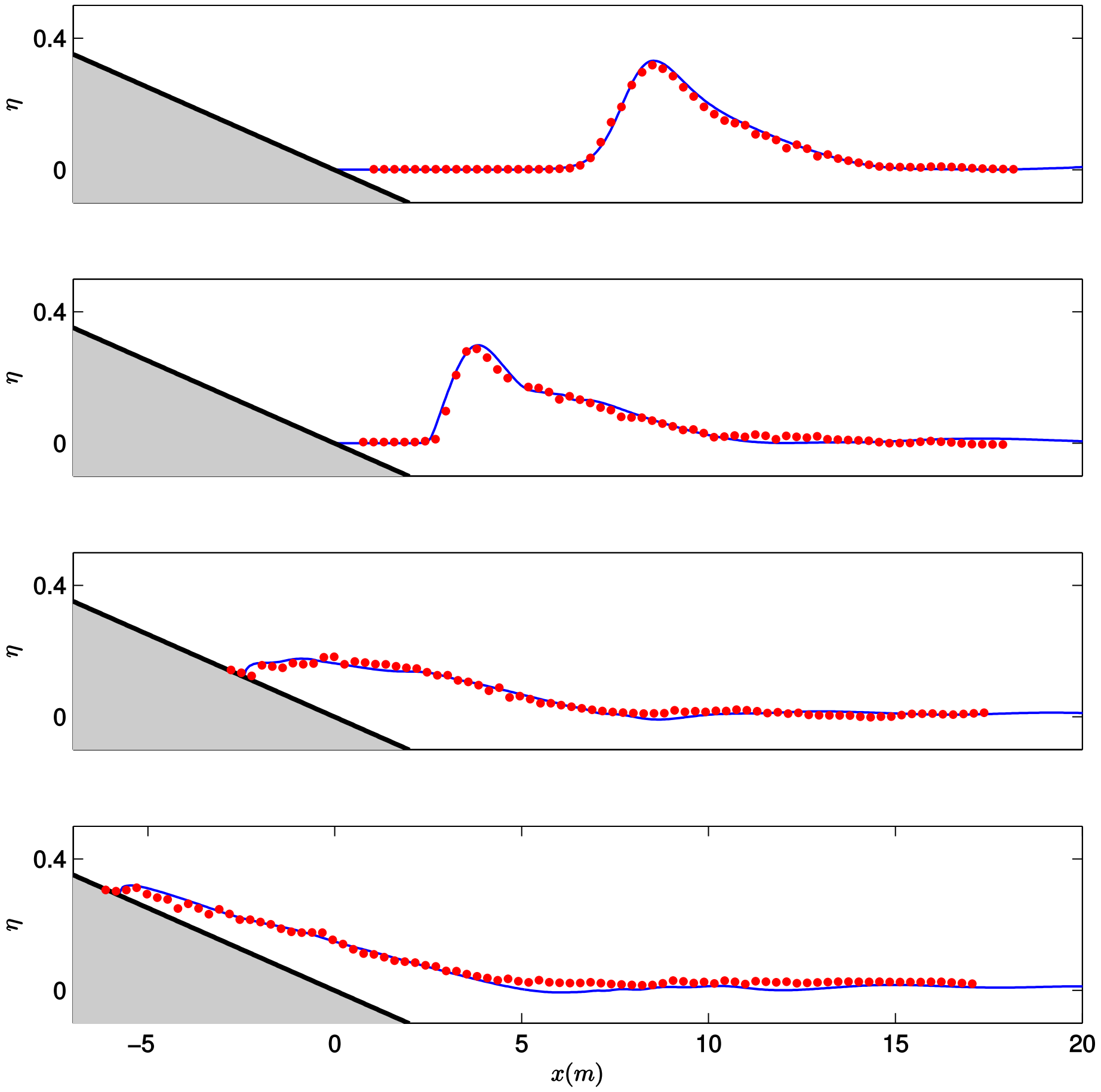

Figure 11: Comparison of experimental data (red) and simulated ones (blue) at times $t\sqrt{g/H} = 15,\ 20,\ 25,\ 30\ s$ with a breaking criteria

The breaking mechanism also works properly in terms of grid convergence. Figure 13 shows a snapshot at time $t\sqrt{g/H} = 15$ for different mesh sizes.

In addition, good results are obtained at maximum run-up, where breaking mechanism also plays an important role. Note that no additional wet-dry treatment on the second step of the scheme is necessary.

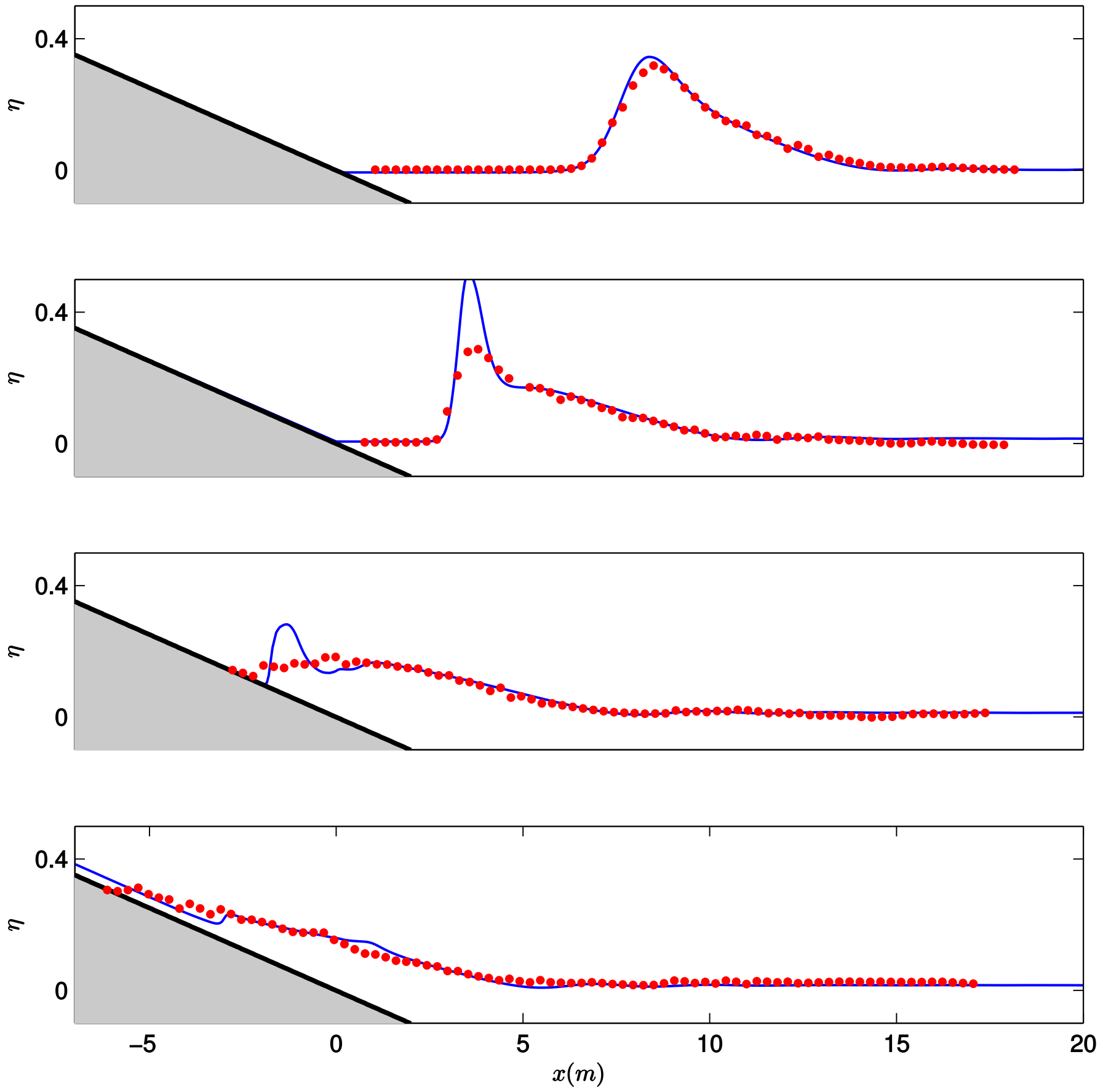

Figure 12: Comparison of experimental data (red) and simulated ones (blue) at times $t\sqrt{g/H} = 15,\ 20,\ 25,\ 30\ s$ without a breaking criteria

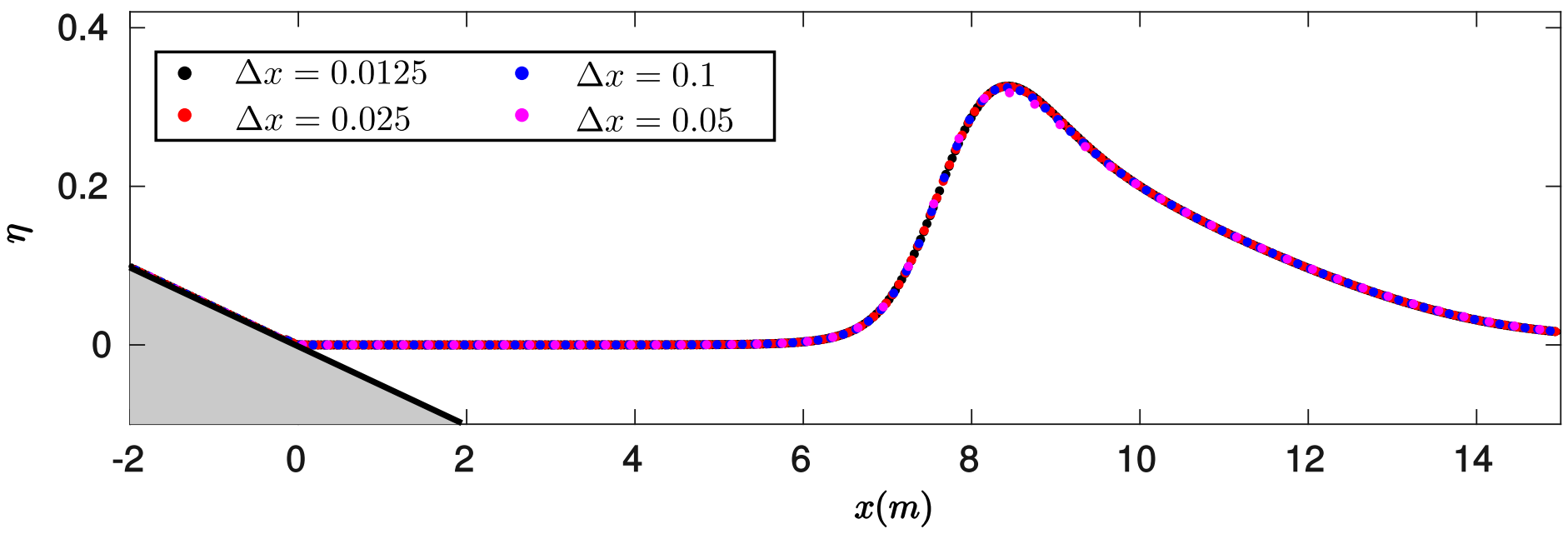

Figure 13: Comparison of free-surface simulation at time $t\sqrt{g/H} = 15$ for different mesh sizes

## 7.5 Solitary wave propagation over reefs

A test case on solitary wave over an idealized fringing reef examines the model's capability of handling nonlinear dispersive waves, breaking waves and bore propagation. The test configurations include a fore reef, a flat reef, and an optional reef crest to represent fringing reefs commonly found in tropical environment. Figure 14 shows a sketch of the laboratory experiments carried out at the O.H. Hinsdale Wave Research Laboratory of Oregon State University. The uni-dimensional domain $[0, 45]$ is discretized with $\Delta x = 0.045\ m$. The bathymetry is defined in the Figure 14.

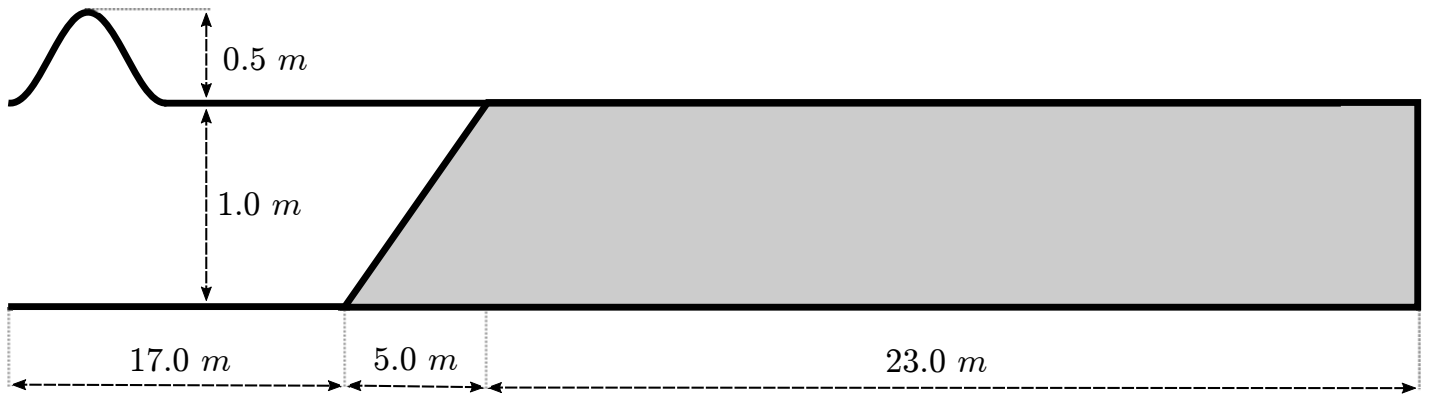

Figure 14: Sketch of the topography

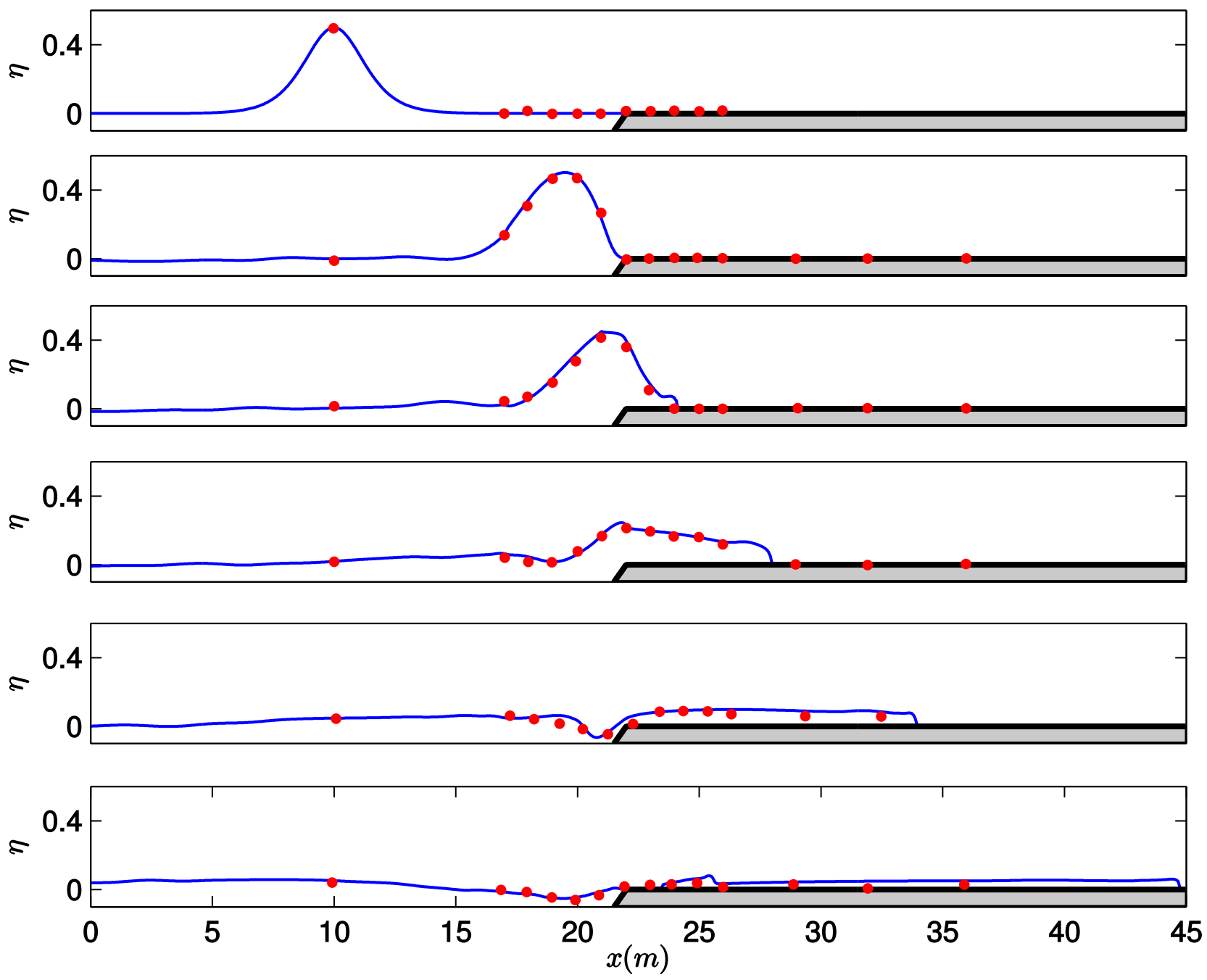

Figure 15: Comparison of experimental data (red points) and numerical (blue) at times $t\sqrt{g/H} = 0,\ 80,\ 100,\ 130,\ 170,\ 250\ s$

A solitary wave of amplitude 0.5 $m$ is placed at point $x = 10\ m$ given by (33). A Manning coefficient of $n_m = 0.012$ was used in order to define the glass surface roughness used in the experiments. Breaking mechanism is considered with $B_1 = 0.15$ and $B_2 = 0.5$. Finally $CFL = 0.9$ and $g = 9.81\ m/s^2$. Free outflow boundary conditions are imposed.

Figure 15 shows snapshots at different times, $t\sqrt{g/H} = t_0$ where $H = 1$. Again, comparison between experimental and simulated data allows us to validate the numerical approach followed here. The water rushes over the flat reef without producing a pronounced bore-shape. The simulation also captures the offshore component of the rarefaction falls, exposing the reef edge, below the initial water level.

## 7.6 Solitary wave on a conical island

The goal of this 2D-numerical test is to compare numerical model results with laboratory measurements. The experiment was carried out at the Coastal and Hydraulic Laboratory, Engineer Research and Development Center of the U.S. Army Corps of Engineers ([6]). The laboratory experiment consists of an idealized representation of Babi Island, in the Flores Sea, in Indonesia. The produced data sets have been frequently used to validate run-up models ([27], [41]).

A directional wave-maker is used to produce planar solitary waves of specified crest lengths and heights. Domain setup consists of a $25 \times 30$ $m$ basin with a conical island situated near the center. The still water level is $H = 0.32$ $m$. The island has a base diameter of 7.2 $m$ , a top diameter of 2.2 $m$, and it is 0.625 $m$ high with a side slope $1:4$. Wave gauges, $\{WG_1,\ WG_2,\ WG_3,\ WG_4\}$, are distributed around the island in order to measure the free surface elevation (see Figure 16).

For the numerical simulation the computational domain is $[-5, 23] \times [0, 28]$ with $\Delta x = 2$ $cm$ and $\Delta y = 2$ $cm$. Free outflow boundary conditions are imposed.

As initial condition for $\eta$ and $u$, a solitary wave (33) of Amplitude $A = 0.06$ $m$ centered at $x = 0$ is given. The wave propagates until 30 $s$, with $CFL = 0.9$ and $g = 9.81$ $m/s^2$. A Manning coefficient of $n_m = 0.015$ is used and breaking mechanism with $B_1 = 0.15$ and $B_2 = 0.5$ is considered.

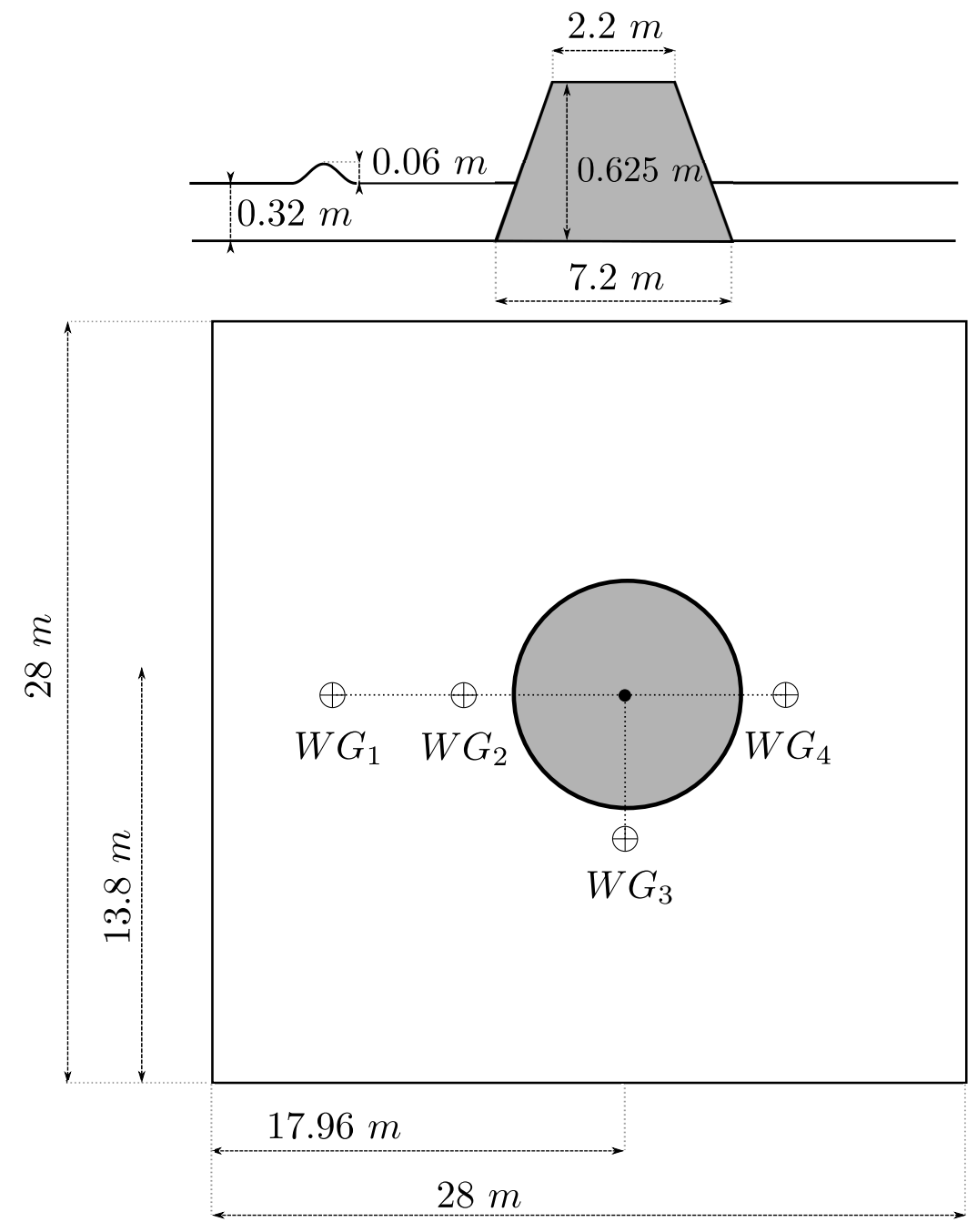

Figure 16: Sketch of the topography

Numerical simulation shows two wave fronts splitting in front of the island and colliding behind it (Figure 19). Comparison with measured and computed water level at gauges $WG_1$, $WG_2$, $WG_3$, $WG_4$ shows good a good agreement. The same is true for the comparison between computed run-up and laboratory measurements (see Figure 17 and Figure 18).

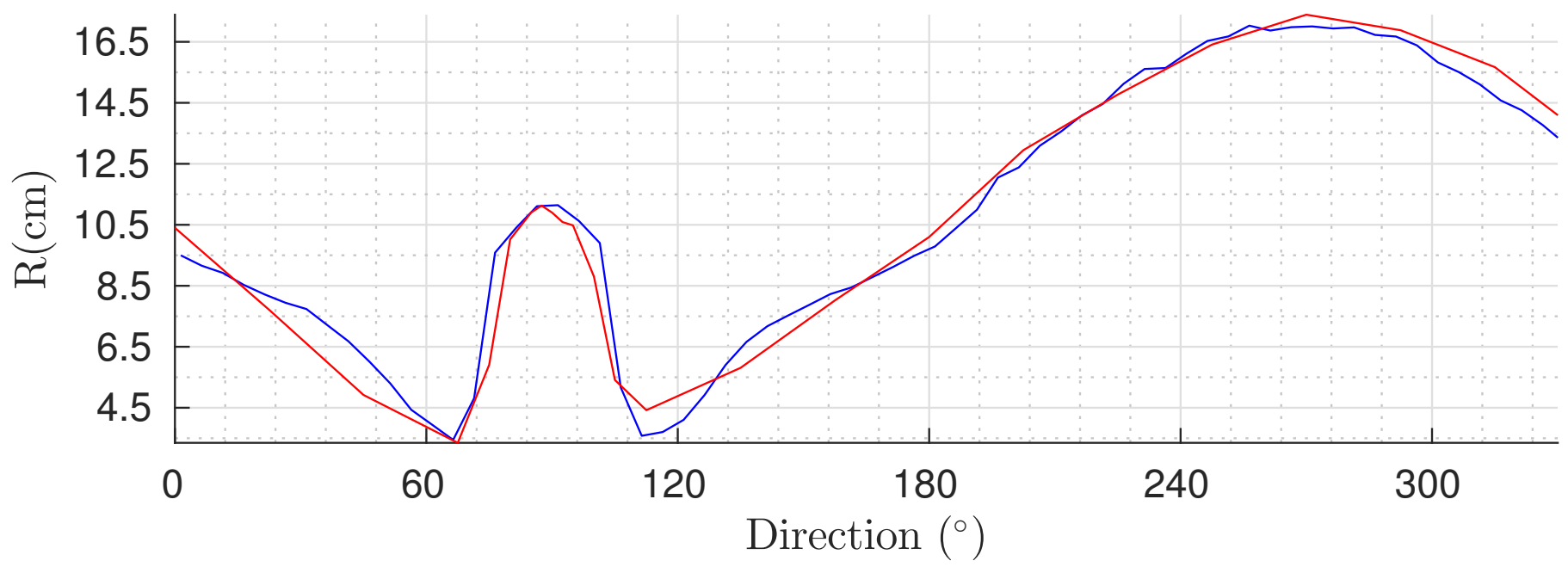

Figure 17: Maximum run-up measured (red) and simulated (blue)

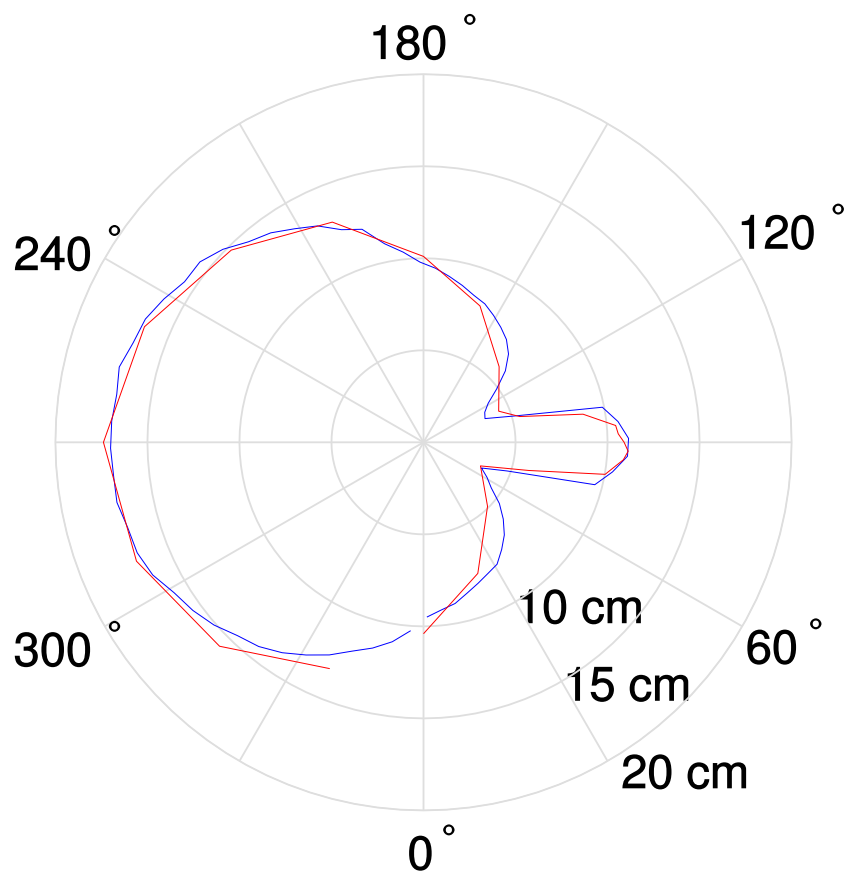

Figure 18: Maximum run-up measured (red) and simulated (blue) in polar coordinates

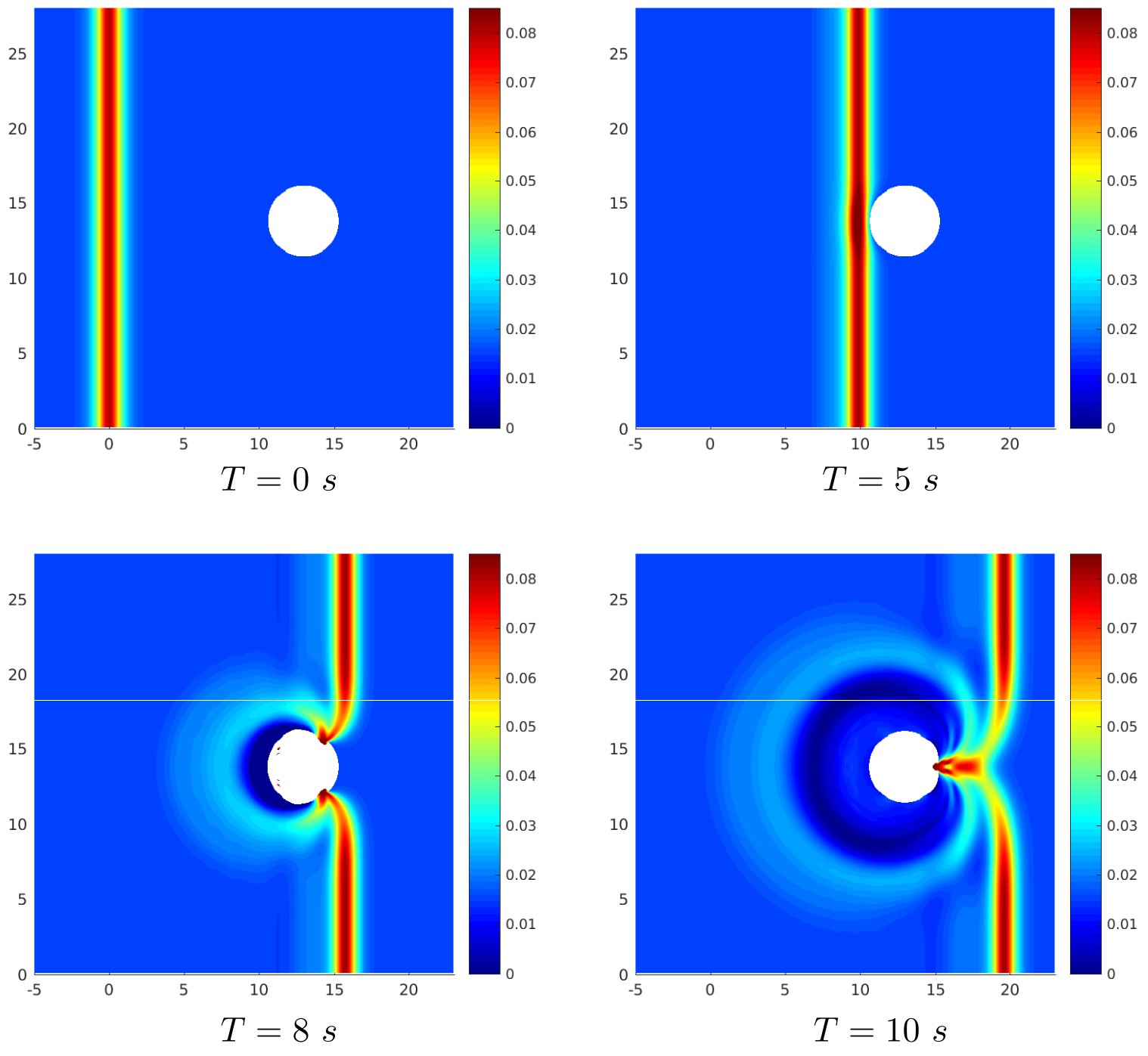

Figure 19: Comparison of numerically calculated free surface $\eta$ at various times.

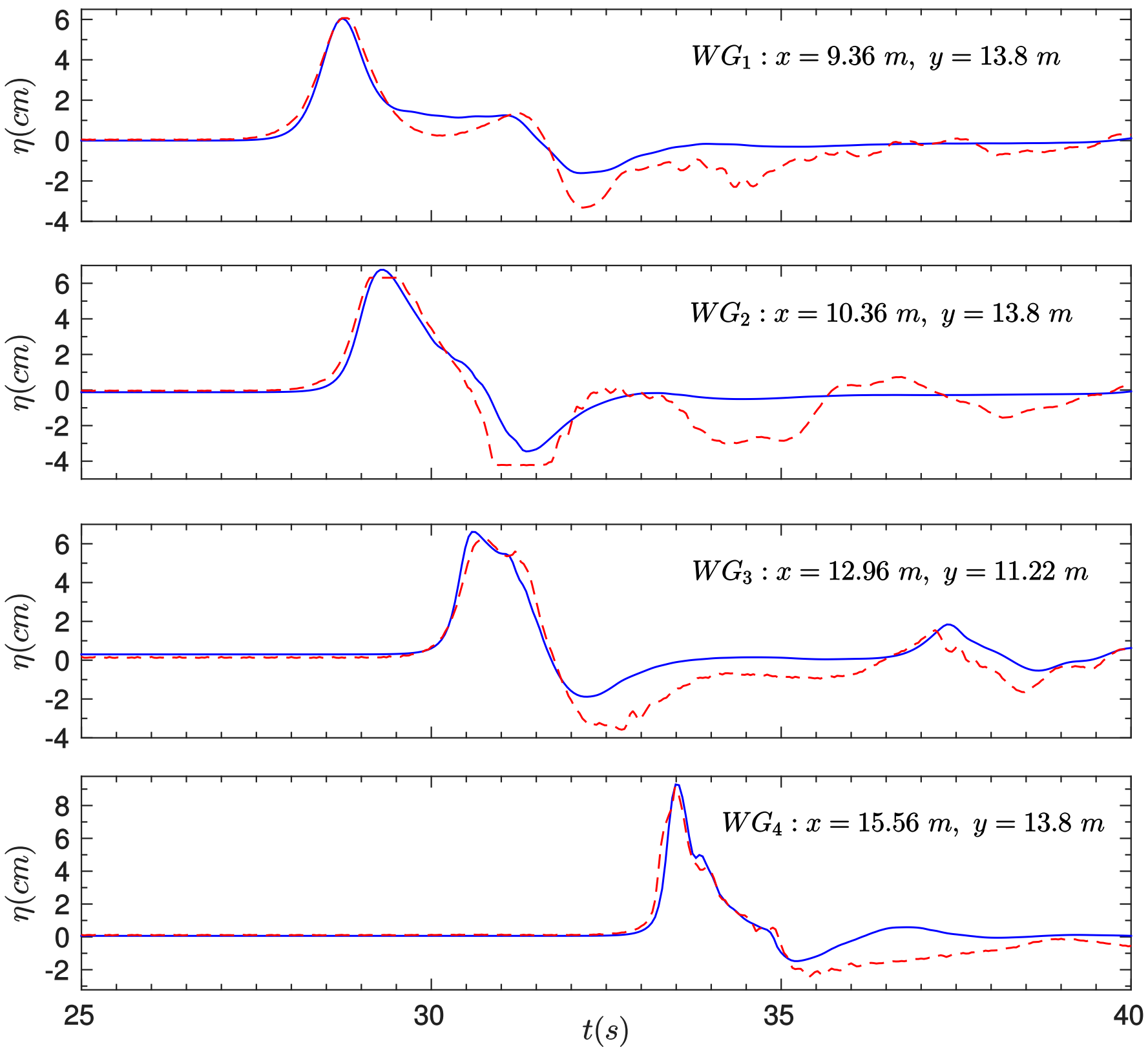

Figure 20: Comparison of data time series (red) and numerical (blue) at wave gauges $WG_1,\ WG_2,\ WG_3,\ WG_4$

## 7.7 Circular dam-break

In this 2D-test case we consider a circular dam-break problem in the $[-5,5]\times[-5,5]$ domain. The depth function is $H(x,y) = 1 - 0.25e^{-x^2-y^2}$ and the initial condition is

$$U_i^0(x,y) = \begin{pmatrix} h^0(x,y) \\ 0 \\ 0 \end{pmatrix}, \quad h^0(x,y) = \begin{cases} H(x,y) & \text{if } \sqrt{x^2+y^2} \leq 0.5 \\ H(x,y) + 0.25 & \text{otherwise.} \end{cases}$$

The goal of this numerical test is to compare the execution times in seconds for the SWE and non-hydrostatic GPU codes for different mesh sizes. Simulations are carried out in the time interval $[0,1]$. The CFL parameter is set to 0.9 and open boundary conditions are considered.

Table 2 shows execution times for both double precision CUDA codes. Different parameters of $\epsilon \in \{10^{-3}, 10^{-4}, 10^{-5}\}$ were taken into account, where $\epsilon$ was defined in (32). Figure (22) shows the results with the different tolerance parameter $\epsilon$. We would like to stress that no big differences are observed for the range of values considered for the tolerance parameter.

Figure 21 shows the speedup achieved using a GPU implementation on a GTX Titan Black with respect to a sequential CPU version of the code. We remark a gain in performance greater than 110.

| Number of Volumes | Runtime (s) | | | |
|---|---|---|---|---|
| | SWE | Non-Hydrostatic | | |
| | | $\epsilon = 10^{-3}$ | $\epsilon = 10^{-4}$ | $\epsilon = 10^{-5}$ |
| $250 \times 250$ | 0.64 | 0.64 | 1.88 | 3.47 |
| $500 \times 500$ | 2.29 | 5.79 | 8.44 | 33.54 |
| $750 \times 750$ | 7.17 | 17.33 | 25.78 | 99.58 |
| $1000 \times 1000$ | 16.75 | 40.47 | 57.23 | 198.91 |
| $1250 \times 1250$ | 33.88 | 79.67 | 143.19 | 381.89 |
| $1500 \times 1500$ | 56.38 | 136.12 | 243.86 | 662.51 |

Table 2: Execution times in sec for SWE and $NH$ GPU implementations

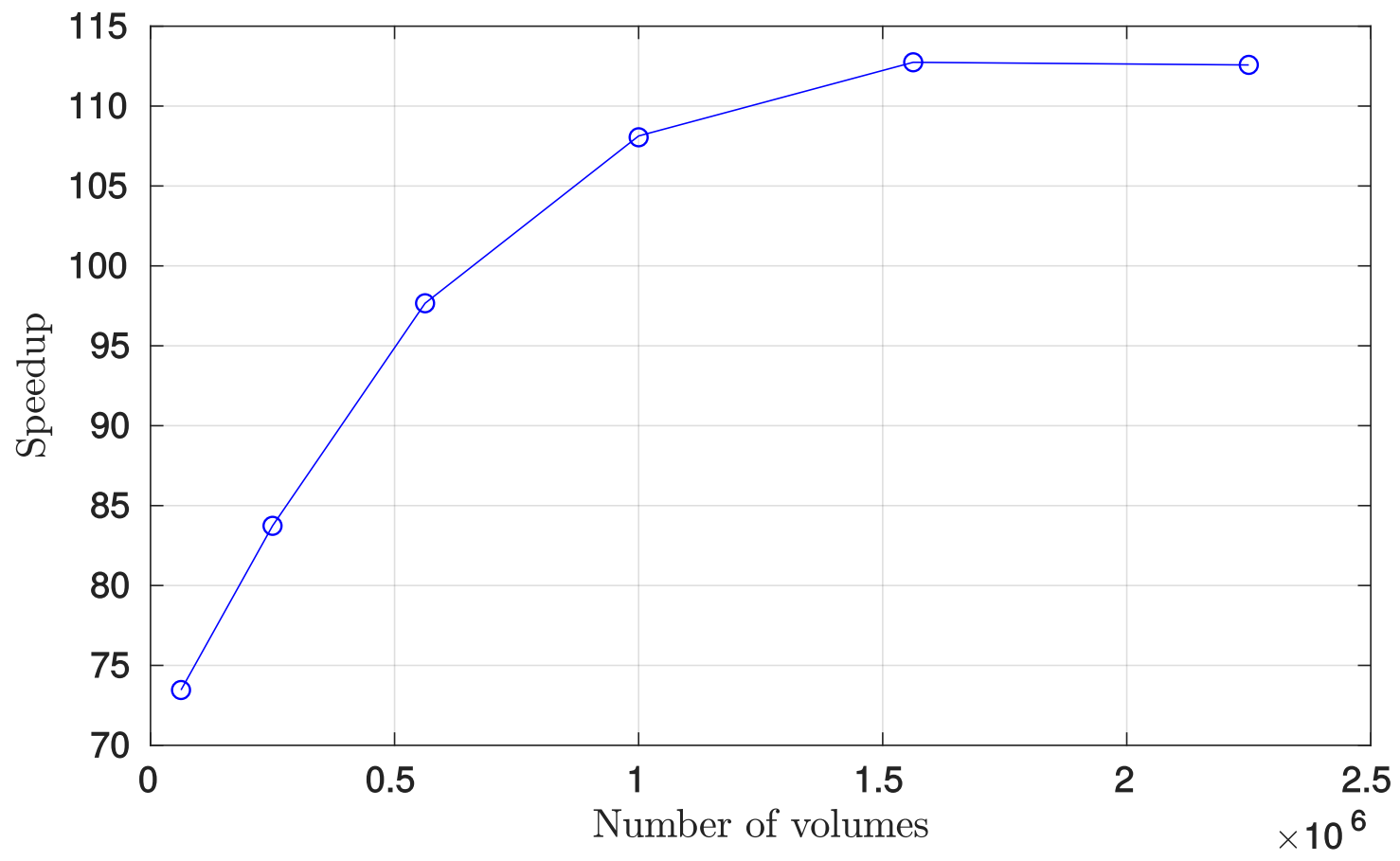

Figure 21: Speedup with respect to a CPU-sequential version of the code

It can be stated thus that the scheme presented here is efficient and can model dispersive effects with a moderate computational cost. To our knowledge, similar models and/or numerical schemes that intend to simulate dispersive effects in such frameworks are much more expensive from the computational point of view.

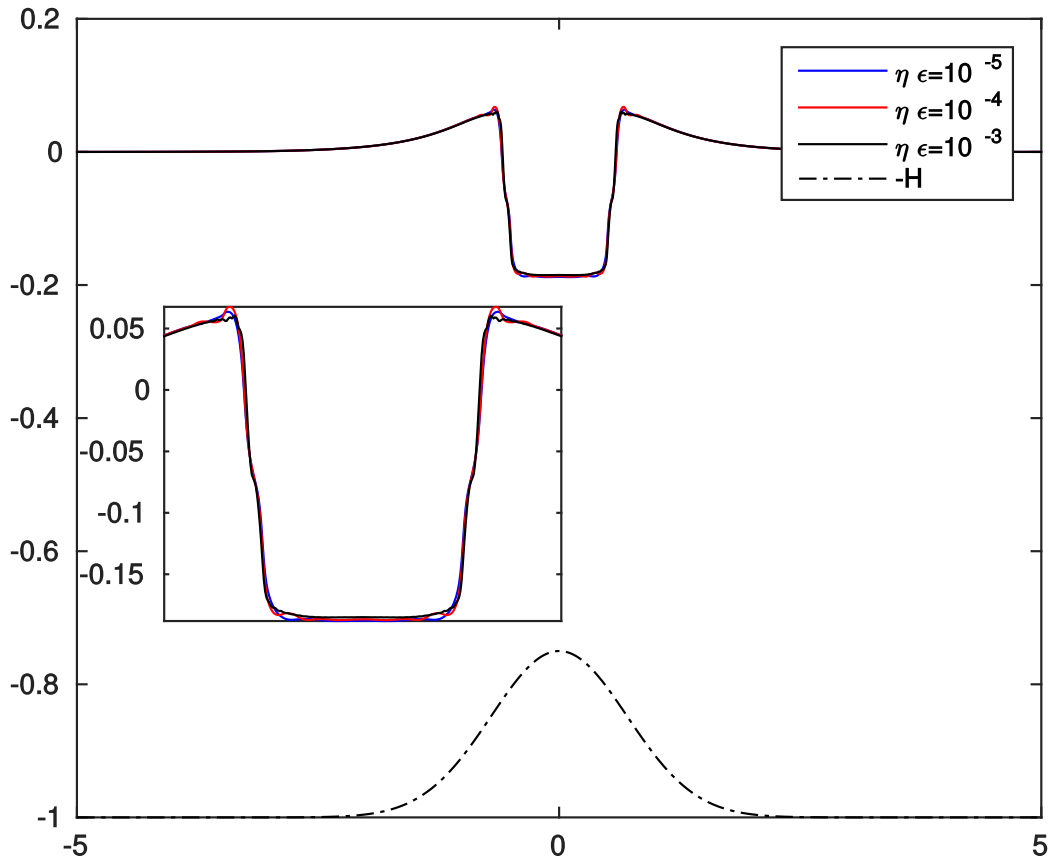

Free surface cross-section at $t = 0.5$ $s$

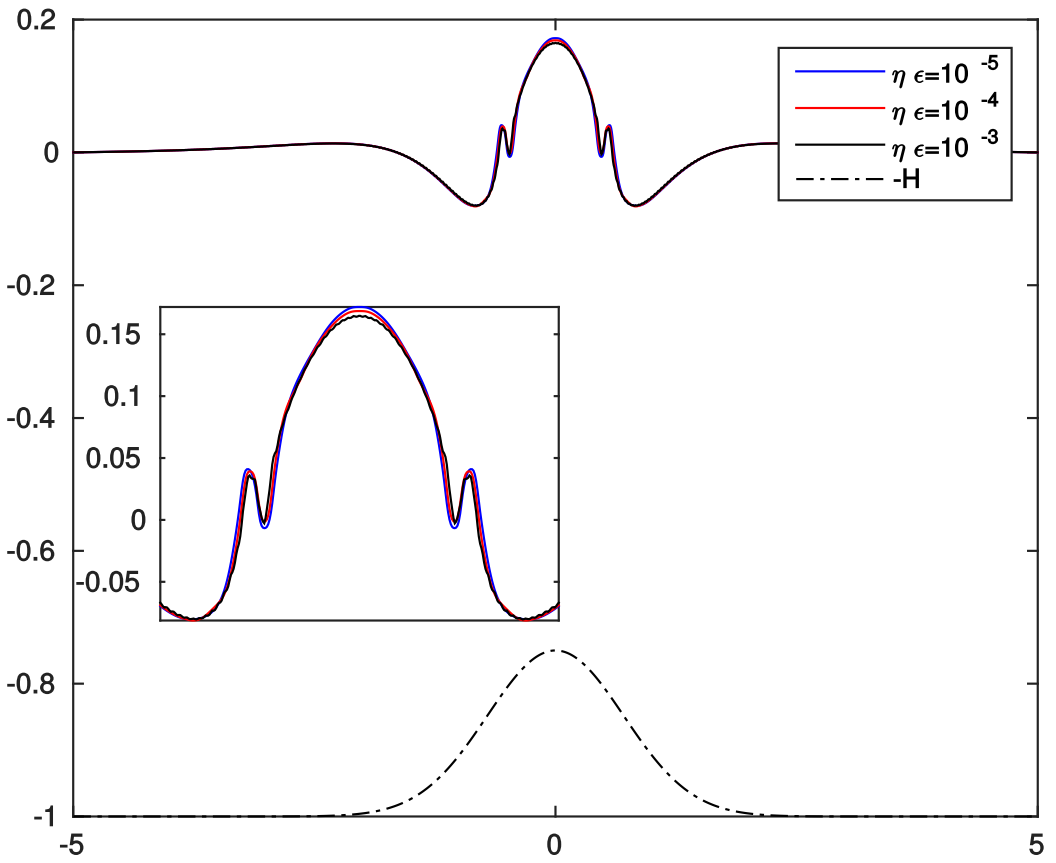

Free surface cross-section at $t = 1.0$ $s$

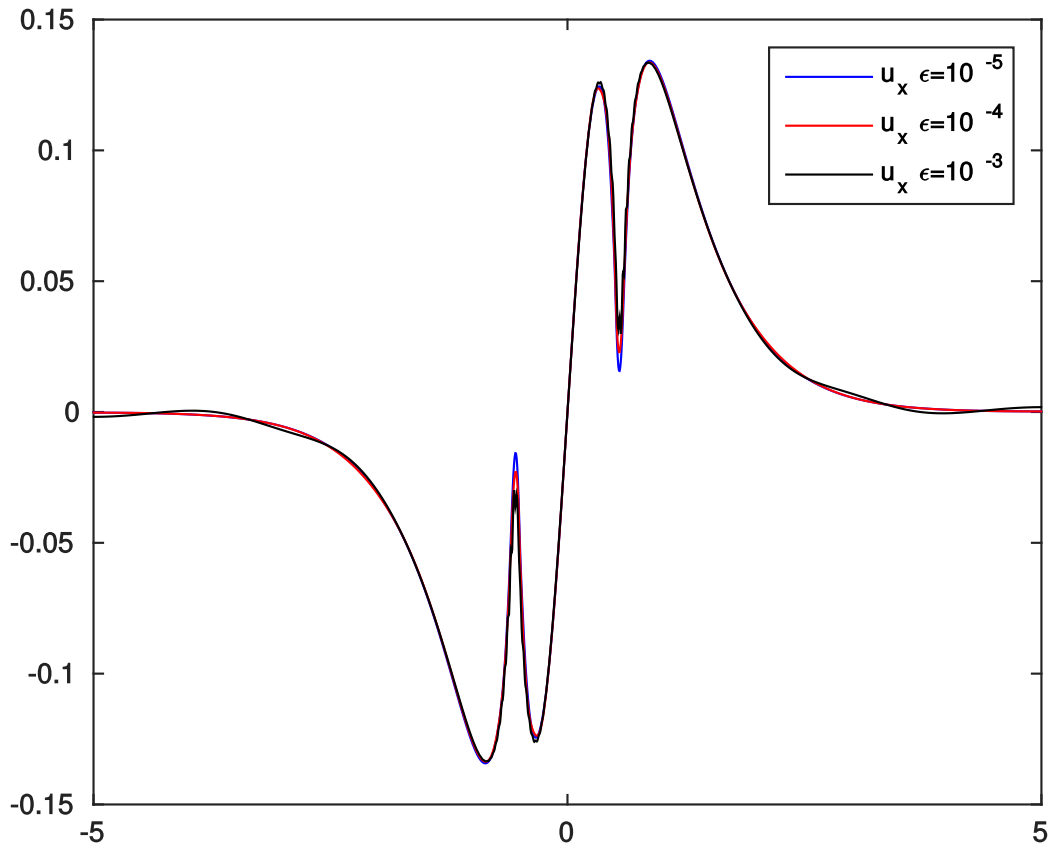

$u_x$ cross-section at $t = 0.5$ $s$

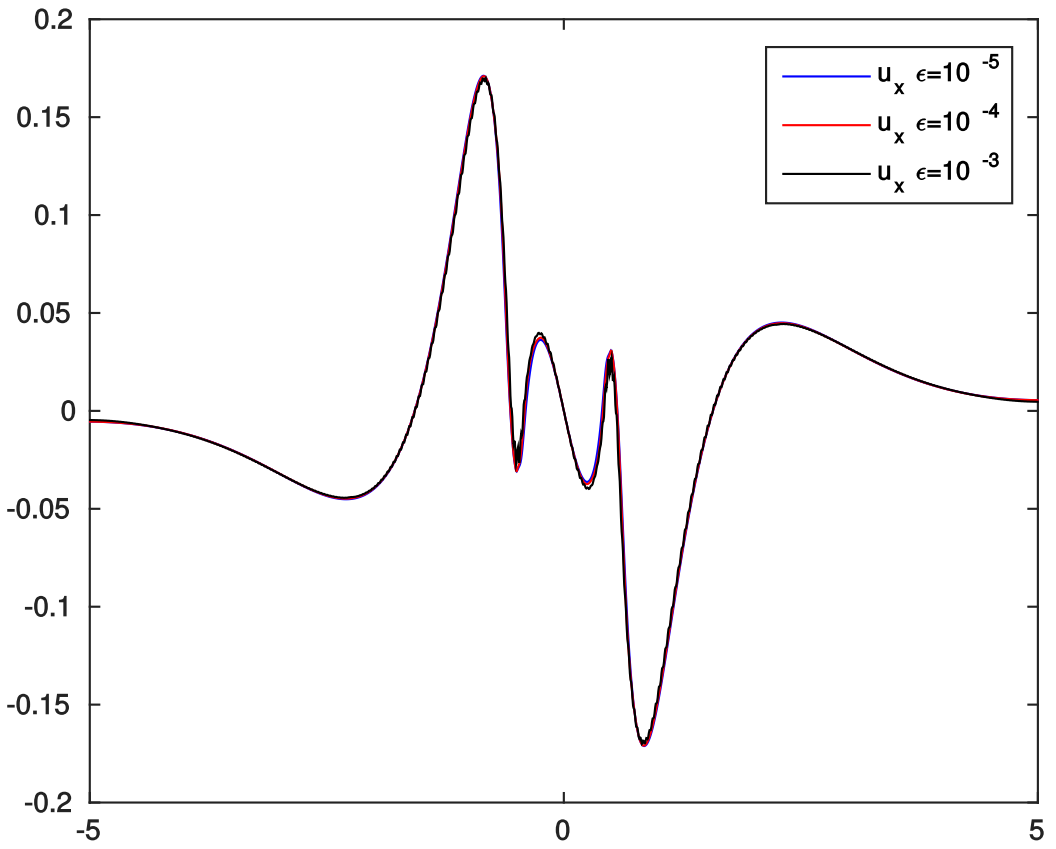

$u_x$ cross-section at $t = 1.0$ $s$

Figure 22: Cross-section of numerical simulations at times $T = 0.5$ $s$ (left) and $T = 1.0$ $s$ (right) for $\epsilon \in \{10^{-3}, 10^{-4}, 10^{-5}\}$.

# 8 Conclusions

In this work, a non-hydrostatic model has been considered in order to incorporate dispersive effects in the propagation of waves in a homogeneous, inviscid and incompressible fluid.

The numerical scheme employed combines a path-conservative finite volume scheme for the underlying hyperbolic system and a finite difference scheme for discretization of non-hydrostatic terms. Furthermore, it is second order accurate and it is well-balanced for the water at rest solution and linearly $L^\infty$-stable under the usual CFL condition.

A wet-dry treatment presented in [9] for the SWE is adopted. Moreover, no numerical truncation for the non-hydrostatic pressure is needed at wet-dry areas, where non-hydrostatic pressure vanishes, as it is usually done (see [41]). This is due to the writing of the equations proposed in (7). To the best of our knowledge, this is an improvement on non-hydrostatic numerical schemes, where usually non-hydrostatic pressure is truncated to zero up to a threshold value.

For such models, it is necessary to consider some dissipative mechanism for breaking waves in order to accurately model waves near the coastal areas. Discretization of the viscosity term needs to solve an extra elliptic problem, which results in additional computational cost. We have proposed a reinterpretation of the viscosity term which results in a new, simple and efficient way to solve the problem. Moreover, the breaking mechanism works adequately in terms of grid-convergence, which is a nice feature as it was exhibited in the numerical test 7.4.

A GPU implementation of the 2D model is carried out. From a computational point of view, the non-hydrostatic code presents good computational times with respect to the SWE GPU times. A numerical test was carried out in order to illustrate such claim. For a tolerance of $\epsilon = 10^{-3}$ for the iterative method that solves the linear system, the wall-clock times for the non-hydrostatic code are no higher than 2.4 times than the SWE code for refined meshes. The achieved speed-up of the GPU implementation, compared with a sequential implementation of the algorithm, is remarkable.

Numerical simulations show that the approach presented here, correctly solves the propagation of solitary waves, preserving their shape for large integration times accurately. Comparison with experimental data is also presented. Experimental data justifies the need to incorporate dispersive effects to faithfully capture waves in the vicinity of the continental shelf. Moreover, complex processes such as run-up, shoaling, wet-dry areas are simulated successfully for the proposed 1D and 2D tests, which validates the approach used here.

The numerical scheme presented in this work provides thus an efficient and accurate approach to model dispersive effects in the propagation of waves near coastal areas.

# A 2D numerical scheme

We consider, as in Section 4, the system:

$$\begin{cases} \boldsymbol{U}_t + \left(\boldsymbol{F}_{1,SW}(\boldsymbol{U})\right)_x + \left(\boldsymbol{F}_{2,SW}(\boldsymbol{U})\right)_y = \\ \boldsymbol{G}_{1,SW}(\boldsymbol{U})H_x + \boldsymbol{G}_{2,SW}(\boldsymbol{U})H_y + \boldsymbol{\mathcal{T}}_{NH}(h, \nabla h, H, \nabla H, p, \nabla p), \\ hw_t = p, \\ \mathcal{B}(\boldsymbol{U}, \nabla h, (\nabla \cdot \boldsymbol{Q}), H, \nabla H, w) = 0, \end{cases} \tag{34}$$

where we denote the vector of the state variables, and the corresponding flows

$$\boldsymbol{U} = \begin{pmatrix} h \\ \boldsymbol{Q} \end{pmatrix}, \quad \boldsymbol{Q} = \begin{pmatrix} q_1 \\ q_2 \end{pmatrix},$$

$$\boldsymbol{F}_{1,SW}(\boldsymbol{U}) = \begin{pmatrix} q_1 \\ \dfrac{q_1^2}{h} + \dfrac{1}{2}gh^2 \\ \dfrac{q_1 q_2}{h} \end{pmatrix}, \quad \boldsymbol{F}_{2,SW}(\boldsymbol{U}) = \begin{pmatrix} q_2 \\ \dfrac{q_1 q_2}{h} \\ \dfrac{q_2^2}{h} + \dfrac{1}{2}gh^2 \end{pmatrix}.$$

The sources terms are given by

$$\boldsymbol{G}_{1,SW}(\boldsymbol{U}) = \begin{pmatrix} 0 \\ gh \\ 0 \end{pmatrix}, \quad \boldsymbol{G}_{2,SW}(\boldsymbol{U}) = \begin{pmatrix} 0 \\ 0 \\ gh \end{pmatrix},$$

and the friction term vector, where a Manning empirical formula is used, is given by

$$\boldsymbol{\tau} = \begin{pmatrix} 0 \\ ghu_1 \dfrac{n^2\sqrt{u_1^2+u_2^2}}{h^{4/3}} \\ ghu_2 \dfrac{n^2\sqrt{u_1^2+u_2^2}}{h^{4/3}} \end{pmatrix}.$$

Finally, non-hydrostatic terms are

$$\boldsymbol{\mathcal{T}}_{NH}(h, \nabla h, H, \nabla H, p, \nabla p) = \begin{pmatrix} 0 \\ \boldsymbol{\mathcal{T}}^{Hor}(h, h_x, H, H_x, p, p_x) \\ \boldsymbol{\mathcal{T}}^{Ver}(h, h_y, H, H_y, p, p_y) \end{pmatrix},$$

being $\boldsymbol{\mathcal{T}}^{Hor}$, $\boldsymbol{\mathcal{T}}^{Ver}$ the horizontal and vertical non-hydrostatic contributions respectively:

$$\boldsymbol{\mathcal{T}}^{Hor}(h, h_x, H, H_x, p, p_x) = -\frac{1}{2}(hp_x + p((2\eta - h)_x)),$$
$$\boldsymbol{\mathcal{T}}^{Ver}(h, h_y, H, H_y, p, p_y) = -\frac{1}{2}(hp_y + p((2\eta - h)_y)),$$

and the free divergence equation is

$$\mathcal{B}(\boldsymbol{U}, \nabla h, (\nabla \cdot \boldsymbol{Q}), H, \nabla H, w) = h\,(\nabla \cdot \boldsymbol{Q}) - \boldsymbol{Q} \cdot \nabla(2\eta - h) + 2hw.$$

We describe now the numerical scheme used to discretize the 2D system (34). The 2D-SWE are written in vector conservative form,

$$\boldsymbol{U}_t + (\boldsymbol{F}_{1,SW}(\boldsymbol{U}))_x + (\boldsymbol{F}_{2,SW}(\boldsymbol{U}))_y = \boldsymbol{G}_{1,SW}(\boldsymbol{U})H_x + \boldsymbol{G}_{2,SW}(\boldsymbol{U})H_y. \quad (35)$$

To discretize (35) the computational domain is decomposed into subsets with a simple geometry, called cells or finite volumes. Here, we consider rectangular structured meshes:

$$V_{ij} = [x_{i-1/2}, x_{i+1/2}] \times [y_{j-1/2}, y_{j+1/2}] \subset \mathbb{R}^2,\ i \in Nx,\ j \in N_y.$$

Given a finite volume $V_{ij}$, $|V_{ij}|$ will represent its area and $\boldsymbol{U}_{ij}(t)$ the constant approximation to the average of the solution in the cell $V_{ij}$ at time $t$ provided by the numerical scheme:

$$\boldsymbol{U}_{ij}(t) = \frac{1}{|V_{ij}|}\int_{V_{ij}} \boldsymbol{U}(\boldsymbol{x}, t)\ d\boldsymbol{x}.$$

Regarding non-hydrostatic terms, we will use one common arrangement of the variables, known as the Arakawa C-grid (see Figure 3). This is an extension of the procedure used for the 1D case. Variables $p$ and $w$ will be computed at the intersections of the edges:

$$p_{i+1/2,j+1/2}(t) = p(x_{i+1/2}, y_{j+1/2}, t),\ w_{i+1/2,j+1/2}(t) = w(x_{i+1/2}, y_{j+1/2}, t),$$

and non-hydrostatic terms will be approximated by second order compact finite differences. The resulting ODE system is discretized using a TVD Runge-Kutta method [18]. For the sake of clarity, only a first order discretization in time will be described. The source terms corresponding to friction terms are discretized semi-implicitly. Thus, friction terms are neglected and only flux, and source terms are considered.

## A.1 Finite volume scheme

For the finite volume scheme we will follow the ideas given in [13] for the two-dimensional problem. In particular, we use the 2D extension of the PVM scheme described in Section 3 in [18].

## A.2 Finite differences scheme

In this subsection, non-hydrostatic variables $p$ and $w$ will be discretized using second order compact finite differences. Following the same procedure as for the 1D equations. Let us define the North and South approximations in the middle of the horizontal edges for the volume $V_{i,j}$ of $\boldsymbol{\mathcal{T}}_{NH}^{Hor}$ by

$$
\begin{aligned}
\boldsymbol{\mathcal{T}}_{N(i,j)}^{Hor}(h, h_x, H, H_x, p, p_x) &= -\frac{1}{2} h_{i,j} \frac{p_{i+1/2,j+1/2} - p_{i-1/2,j+1/2}}{\Delta x} \\
&\quad - \frac{1}{2} \frac{p_{i+1/2,j+1/2} + p_{i-1/2,j+1/2}}{2} \cdot \frac{2\eta_{i+1,j} - h_{i+1,j} - (2\eta_{i-1,j} - h_{i-1,j})}{2\Delta x}, \\
\boldsymbol{\mathcal{T}}_{S(i,j)}^{Hor}(h, h_x, H, H_x, p, p_x) &= -\frac{1}{2} h_{i,j} \frac{p_{i+1/2,j-1/2} - p_{i-1/2,j-1/2}}{\Delta x} \\
&\quad - \frac{1}{2} \frac{p_{i+1/2,j-1/2} + p_{i-1/2,j-1/2}}{2} \cdot \frac{2\eta_{i+1,j} - h_{i+1,j} - (2\eta_{i-1,j} - h_{i-1,j})}{2\Delta x},
\end{aligned}
$$

respectively.

Same ideas for the East and West approximations in the middle of the vertical edges for the volume $V_{i,j}$ of $\boldsymbol{\mathcal{T}}_{NH}^{Ver}$:

$$
\begin{aligned}
\boldsymbol{\mathcal{T}}_{E(i,j)}^{Ver}(h, h_y, H, H_y, p, p_y) &= -\frac{1}{2} h_{i,j} \frac{p_{i+1/2,j+1/2} - p_{i+1/2,j-1/2}}{\Delta y} \\
&\quad - \frac{1}{2} \frac{p_{i+1/2,j+1/2} + p_{i+1/2,j-1/2}}{2} \cdot \frac{2\eta_{i,j+1} - h_{i,j+1} - (2\eta_{i,j-1} - h_{i,j-1})}{2\Delta y}, \\
\boldsymbol{\mathcal{T}}_{W(i,j)}^{Ver}(h, h_y, H, H_y, p, p_y) &= -\frac{1}{2} h_{i,j} \frac{p_{i-1/2,j+1/2} - p_{i-1/2,j-1/2}}{\Delta y} \\
&\quad - \frac{1}{2} \frac{p_{i-1/2,j+1/2} + p_{i-1/2,j-1/2}}{2} \cdot \frac{2\eta_{i,j+1} - h_{i,j+1} - (2\eta_{i,j-1} - h_{i,j-1})}{2\Delta y}.
\end{aligned}
$$

Note that, if we approximate

$$
\boldsymbol{\mathcal{T}}_{NH}(h, \nabla(h), H, \nabla(H), p, \nabla(p))_{i,j} \approx \begin{pmatrix} 0 \\ \frac{1}{2}\left(\boldsymbol{\mathcal{T}}_{N(i,j)}^{Hor} + \boldsymbol{\mathcal{T}}_{S(i,j)}^{Hor}\right) \\ \frac{1}{2}\left(\boldsymbol{\mathcal{T}}_{E(i,j)}^{Ver} + \boldsymbol{\mathcal{T}}_{W(i,j)}^{Ver}\right) \end{pmatrix}, \qquad (36)
$$

then we have a second order approximation of $\boldsymbol{\mathcal{T}}_{NH}(h, \nabla(h), H, \nabla(H), p, \nabla(p))$ at the center of the volume $V_{i,j}$.

Likewise, $\mathcal{B}(\boldsymbol{U}, \nabla(h), (\nabla\cdot\boldsymbol{Q}), H, \nabla(H), w)$ will be discretized for every point $(x_{j+1/2}, y_{i+1/2})$ of the staggered mesh by

$$\mathcal{B}(\boldsymbol{U}, \nabla(h), (\nabla\cdot\boldsymbol{Q}), H, \nabla(H), w)_{i+1/2,j+1/2} \approx h_{i+1/2,j+1/2}(\nabla\cdot\boldsymbol{Q})_{i+1/2,j+1/2} - \boldsymbol{Q}_{i+1/2,j+1/2}\cdot\nabla(2\eta - h)_{i+1/2,j+1/2} + 2h_{i+1/2,j+1/2}w_{i+1/2,j+1/2}, \tag{37}$$

being

$$h_{i+1/2,j+1/2} = \frac{1}{4}\left(h_{i,j} + h_{i+1,j} + h_{i+1,j+1} + h_{i,j+1}\right), \tag{38}$$

$$(\nabla\cdot\boldsymbol{Q})_{i+1/2,j+1/2} = \frac{q_{1,E} - q_{1,W}}{\Delta x} + \frac{q_{2,N} - q_{2,S}}{\Delta y},$$

$$\boldsymbol{Q}_{i+1/2,j+1/2} = \begin{pmatrix} \dfrac{q_{1,E} + q_{1,W}}{2} \\ \dfrac{q_{2,N} + q_{2,S}}{2} \end{pmatrix}, \tag{39}$$

$$\nabla(2\eta - h)_{i+1/2,j+1/2} = \begin{pmatrix} \dfrac{(2\eta - h)_E - (2\eta - h)_W}{2} \\ \dfrac{(2\eta - h)_N - (2\eta - h)_S}{2} \end{pmatrix}, \tag{40}$$

where $q_{1,E},\ q_{1,W},\ q_{2,N},\ q_{2,S}$ and $(2\eta - h)_E,\ (2\eta - h)_W,\ (2\eta - h)_N,\ (2\eta - h)_S$ are second order approximations of $q_1,\ q_2$ and $(2\eta - h)$ respectively in the middle of the edges (see Figure(3)). Expressions for this approximations will be introduced in the next section.

### Final Numerical Scheme

Let be given time steps $\Delta t^n$, note $t^n = \sum_{k\le n}\Delta t^k$ and $\mathbf{U}_{i,j}(t^n) = \mathbf{U}^n_{i,j}$, $p_{i+1/2}(t^n) = p^n_{i+1/2}$, $w_{i+1/2}(t^n) = w^n_{i+1/2}$. The proposed numerical scheme consists of two steps:

On a first stage, SWE approximation is carried out. Let us define $U_{i,j}^{n+1/2}$ as the averaged values of $U$ on cell $I_i$ at time $t^n$ for the SWE as detailed in the subsection (A.1).

On a second stage, we consider the system

$$\begin{cases} \boldsymbol{U}_{i,j}^{n+1} = \boldsymbol{U}_{i,j}^{n+1/2} + \Delta t \boldsymbol{\mathcal{T}}_{NH}\left(h^{n+1}, \nabla h^{n+1}, H, \nabla H, p^{n+1}, \nabla p^{n+1}\right)_{i,j} \\ w_{i+1/2,j+1/2}^{n+1} = w_{i+1/2,j+1/2}^{n} + \Delta t \dfrac{p_{i+1/2,j+1/2}^{n+1}}{h_{i+1/2,j+1/2}^{n+1}} \\ \mathcal{B}\left(\widetilde{\boldsymbol{U}}^{n+1}, \nabla h^{n+1}, (\nabla \cdot \widetilde{\boldsymbol{Q}}^{n+1}), H, \nabla H, w^{n+1}\right)_{i+1/2,j+1/2} = 0, \end{cases} \tag{41}$$

where:

$$\boldsymbol{\mathcal{T}}_{NH}(h^{n+1}, \nabla h^{n+1}, H, \nabla H, p^{n+1}, \nabla p^{n+1})_{i,j}$$

is defined by (36),

$$h_{i+1/2,j+1/2}$$

is defined by (38) and

$$\mathcal{B}(\widetilde{\boldsymbol{U}}^{n+1}, \nabla h^{n+1}, (\nabla \cdot \widetilde{\boldsymbol{Q}}^{n+1}), H, \nabla H, w^{n+1})_{i+1/2,j+1/2}$$

is defined by (37), being

$$\begin{aligned} q_{1,E}^{n+1} &= \frac{1}{2}\left(q_{x,i+1,j+1}^{n+1} + q_{x,i+1,j}^{n+1}\right) \\ &+ \frac{1}{2}\Delta t \boldsymbol{\mathcal{T}}_{S(i+1,j+1)}^{Hor}(h^{n+1}, h_y^{n+1}, H, H_y, p^{n+1}, p_y^{n+1}) \\ &+ \frac{1}{2}\Delta t \boldsymbol{\mathcal{T}}_{N(i+1,j)}^{Hor}(h^{n+1}, h_y^{n+1}, H, H_y, p^{n+1}, p_y^{n+1}), \end{aligned}$$

$$\begin{aligned} q_{1,W}^{n+1} &= \frac{1}{2}\left(q_{x,i,j+1}^{n+1} + q_{x,i,j}^{n+1}\right) \\ &+ \frac{1}{2}\Delta t \boldsymbol{\mathcal{T}}_{S(i,j+1)}^{Hor}(h^{n+1}, h_y^{n+1}, H, H_y, p^{n+1}, p_y^{n+1}) \\ &+ \frac{1}{2}\Delta t \boldsymbol{\mathcal{T}}_{N(i,j)}^{Hor}(h^{n+1}, h_y^{n+1}, H, H_y, p^{n+1}, p_y^{n+1}), \end{aligned}$$

$$q_{2,N}^{n+1} = \frac{1}{2}\left(q_{y,i+1,j+1}^{n+1} + q_{x,i,j+1}^{n+1}\right)$$
$$+ \frac{1}{2}\Delta t \boldsymbol{\mathcal{T}}_{W(i+1,j+1)}^{Ver}(h^{n+1}, h_y^{n+1}, H, H_y, p^{n+1}, p_y^{n+1})$$
$$+ \frac{1}{2}\Delta t \boldsymbol{\mathcal{T}}_{E(i,j+1)}^{Ver}(h^{n+1}, h_y^{n+1}, H, H_y, p^{n+1}, p_y^{n+1}),$$

$$q_{2,S}^{n+1} = \frac{1}{2}\left(q_{y,i+1,j}^{n+1} + q_{x,i,j}^{n+1}\right)$$
$$+ \frac{1}{2}\Delta t \boldsymbol{\mathcal{T}}_{W(i+1,j)}^{Ver}(h^{n+1}, h_y^{n+1}, H, H_y, p^{n+1}, p_y^{n+1})$$
$$+ \frac{1}{2}\Delta t \boldsymbol{\mathcal{T}}_{E(i,j)}^{Ver}(h^{n+1}, h_y^{n+1}, H, H_y, p^{n+1}, p_y^{n+1}),$$

$$(2\eta - h)_E^{n+1} = \frac{2\eta_{i+1,j+1} - h_{i+1,j+1} + (2\eta_{i+1,j} - h_{i+1,j})}{2},$$

$$(2\eta - h)_W^{n+1} = \frac{2\eta_{i,j+1} - h_{i,j+1} + (2\eta_{i,j} - h_{i,j})}{2},$$

$$(2\eta - h)_N^{n+1} = \frac{2\eta_{i+1,j+1} - h_{i+1,j+1} + (2\eta_{i,j+1} - h_{i,j+1})}{2},$$

$$(2\eta - h)_S^{n+1} = \frac{2\eta_{i+1,j} - h_{i+1,j} + (2\eta_{i,j} - h_{i,j})}{2}.$$

# B Coefficients and matrix of the linear system

## A.1 Coefficients for the one-dimensional case

The linear system defined in (29)

$$\boldsymbol{A}^{n+1/2}\boldsymbol{\mathcal{P}}^{n+1} = \boldsymbol{\mathcal{RHS}}^{n+1/2},$$

where

$$\boldsymbol{\mathcal{P}}^{n+1} = \begin{pmatrix} p_{1/2}^{n+1} \\ p_{1+1/2}^{n+1} \\ \vdots \\ p_{N+1/2}^{n+1} \end{pmatrix}$$

is given by:

$$\boldsymbol{A}^{n+1/2} = \begin{pmatrix} b_0^{n+1/2} & c_0^{n+1/2} & & & & \cdots & 0 \\ a_1^{n+1/2} & b_1^{n+1/2} & c_1^{n+1/2} & & & & \\ & \ddots & \ddots & \ddots & & & \vdots \\ & & a_i^{n+1/2} & b_i^{n+1/2} & c_i^{n+1/2} & & \\ \vdots & & & \ddots & \ddots & \ddots & \\ & & & & a_{N-1}^{n+1/2} & b_{N-1}^{n+1/2} & c_{N-1}^{n+1/2} \\ 0 & & \cdots & & & a_N^{n+1/2} & b_N^{n+1/2} \end{pmatrix}, \quad (42)$$

where for $k \in \{0, \ldots, N\}$, neglecting the dependence on time in the notation:

$$\begin{cases} a_i = (\xi_{\Delta x,i} - 2h_i)(\xi_{\Delta x,i+1/2} + 2h_{i+1/2}), \\ b_i = 16\Delta x^2 + \xi_{\Delta x,i+1/2}(\xi_{\Delta x,i} + \xi_{\Delta x,i+1} + 2h_i - 2h_{i+1}) + 2h_{i+1/2}(\xi_{\Delta x,i} - \xi_{\Delta x,i+1} + 4h_{i+1/2}), \\ c_i = (\xi_{\Delta x,i+1} + 2h_{i+1})(\xi_{\Delta x,i+1/2} - 2h_{i+1/2}). \end{cases} \quad (43)$$

The coefficients described above are conveniently modified depending on the choice of the boundary conditions.

$h_{i+1/2}$ is given by (25) and

$$\xi_{\Delta x,i} = \Delta x \left(2\eta_{x,i} - h_{x,i}\right), \ \xi_{\Delta x,i+1/2} = \Delta x \left(2\eta_{x,i+1/2} - h_{x,i+1/2}\right),$$

being $\eta_{x,i}$ and $h_{x,i}$ given by (27) and $\eta_{x,i+1/2}$ and $h_{x,i+1/2}$ given by (25).

Finally, the *Right Hand Side* is given by

$$(\boldsymbol{\mathcal{RHS}})_i = \frac{8\Delta x^2}{\Delta t}\left(h_{i+1/2}q_{x,i+1/2} - q_{i+1/2}\left(2\eta_{x,i+1/2} - h_{x,i+1/2}\right) + 2h_{i+1/2}w_{i+1/2}\right),$$

where $q_{i+1/2}$ and $q_{x,i+1/2}$ are given by (25).

## A.2 Analysis of the linear system for small water heights

If we assume

$$h = \delta,\ q = w = 0,\ H = \alpha x$$

then the coefficients (43) reduce to

$$\begin{cases} a_i = 4(\alpha - \delta)(\alpha + \delta), \\ b_i = 8(2\Delta x^2 + \alpha^2 + \delta^2), \\ c_i = 4(\alpha - \delta)(\alpha + \delta), \end{cases}$$

and the *Right Hand Side* vector vanishes

$$\boldsymbol{\mathcal{RHS}} = \mathbf{0}$$

Moreover, since the linear system is strictly diagonal dominant, the matrix $\boldsymbol{A}$ is invertible.

# Acknowledgements

This research has been supported by the Spanish Government through the Research projects MTM2015-70490-C2-1-R, MTM2015-70490-C2-2-R.